\documentclass[12pt,centertags,oneside]{amsart}

\usepackage{amsmath,amstext,amsthm,amscd,typearea,hyperref}
\usepackage{amssymb}
\usepackage{a4wide}
\usepackage[mathscr]{eucal}
\usepackage{mathrsfs}
\usepackage{typearea}
\usepackage{charter}
\usepackage{pdfsync}
\usepackage[a4paper,width=16.2cm,top=3cm,bottom=3cm]{geometry}

\numberwithin{equation}{section}

\newtheorem{theorem}{Theorem}[section]
\newtheorem{definition}[theorem]{Definition}
\newtheorem{proposition}[theorem]{Proposition}
\newtheorem{corollary}[theorem]{Corollary}
\newtheorem{lemma}[theorem]{Lemma}
\newtheorem{remark}[theorem]{Remark}

\newtheorem{example}[theorem]{Example}

\newcommand{\cali}[1]{\mathscr{#1}}

\newcommand{\Tan}{\mathop{\mathrm{Tan}}\nolimits}

\newcommand{\supp}{{\rm supp}}
\newcommand{\const}{{\rm const}}
\newcommand{\dist}{{\rm dist}}

\newcommand{\ddc}{{dd^c}}
\newcommand{\dc}{{d^c}}
\newcommand{\dbar}{{\overline\partial}}
\newcommand{\ddbar}{{\partial\overline\partial}}

\newcommand{\id}{{\rm id}}

\newcommand{\Ac}{\cali{A}}
\newcommand{\Bc}{\cali{B}}
\newcommand{\Cc}{\cali{C}}

\newcommand{\Fc}{\cali{F}}

\newcommand{\FS}{{\rm FS}}

\newcommand{\A}{\mathbb{A}}
\newcommand{\C}{\mathbb{C}}
\newcommand{\D}{\mathbb{D}}
\renewcommand{\H}{\mathbb{H}}

\newcommand{\Z}{\mathbb{Z}}
\newcommand{\R}{\mathbb{R}}

\newcommand{\B}{\mathbb{B}}
\newcommand{\U}{\mathbb{U}}
\renewcommand{\S}{\mathbb{S}}
\renewcommand{\P}{\mathbb{P}}

\renewcommand{\L}{{\mathcal L}}

\newcommand\raisepunct[1]{\,\mathpunct{\raisebox{0.5ex}{#1}}}

\title{Unique ergodicity for foliations in $\P^2$  with an invariant curve}

\author{Tien-Cuong Dinh}
\address{Department of Mathematics, National University 
of Singapore, 10 Lower Kent Ridge Road, Singapore 119076.}
\email{matdtc@nus.edu.sg}
\thanks{T.-C.\ D.\  partially supported by Start-Up 
Grant R-146-000-204-133 from National University of\break Singapore}

\author{Nessim Sibony}
\address{Laboratoire de Math\'ematiques d'Orsay, Univ. Paris-Sud, CNRS, Universit\'e
Paris-Saclay, 91405 Orsay, France.}
\email{Nessim.Sibony@math.u-psud.fr}

\date{September 25, 2015}

\begin{document}

\maketitle

{\it To Duong Hong Phong for his sixtieth birthday}

\begin{abstract}
Consider a foliation in the projective plane admitting a projective line as the unique invariant algebraic curve. Assume that the foliation is generic in the sense that its singular points are hyperbolic.  We show that there is a unique positive $\ddc$-closed $(1,1)$-current of mass 1 which is directed by the foliation and this is the current of integration on the invariant line.  A unique ergodicity theorem for the distribution of leaves follows: for any leaf $L$, appropriate averages of $L$ converge to the current of integration on the invariant line.
The result uses an extension of our theory of densities for currents. Foliations on compact K\"ahler surfaces with one or more invariant curves are also considered.
\end{abstract}

\medskip\medskip

\noindent
{\bf Classification AMS 2010:} 37F75, 37A.

 \medskip


\section{Introduction} \label{introduction}

Let $\Bc_d$ denote the space of holomorphic foliations of degree $d\geq 2$ in $\P^2$. Recall that {\it the (geometric) degree} of a foliation $\Fc$ is the number of tangency points of leaves of $\Fc$ with a generic projective line in $\P^2$. 
A theorem by Jouanolou says that generic foliations in $\Bc_d$ have no leaf contained in an algebraic curve \cite{Jouanolou}. 
On the other hand, let $\Ac_d$ be the space of foliations in $\P^2$ defined in an affine chart $\C^2$ by a polynomial 1-form
$\Omega=Pdz_1-Qdz_2$, with $P,Q$ relatively prime polynomials in $(z_1,z_2)\in\C^2$  such that $\max(\deg P,\deg Q)=d$. The number $d$ is called {\it the affine degree}. The class $\Ac_d$ has been intensely studied by Ilyashenko and Khudai-Veronov among many others. The book by Ilyashenko-Yakovenko \cite{IY} is a good reference for results on that class, see also \cite{I} and Shcherbakov \cite{Sh}. 

It is not difficult to check that 
$$\Ac_d\subset \Bc_{d-1}\cup \Bc_d \quad \text{and} \quad \Bc_d\subset \Ac_d\cup \Ac_{d+1}.$$
Moreover, $\Ac_d$ and $\Bc_d$ can be identified with  Zariski dense open sets of some projective spaces. We will say that a property is {\it  typical} for $\Ac_d$ (resp. $\Bc_d$) if it is valid for $\Fc$ in a set of full Lebesgue measure in $\Ac_d$ (resp. $\Bc_d$). 
We refer to \cite{IY} for the following typical properties of  foliations $\Fc\in \Ac_d$.

\begin{enumerate}
\item The line at infinity $L_\infty$ of $\P^2$ is invariant for $\Fc$; 
this is the case when $z_1P-z_2Q$ is of degree $d+1$.
\item The line at infinity contains $d+1$ distinct singular points (Petrovski-Landis).
\item The line at infinity is the only algebraic variety containing a leaf.
\item All leaves except $L_\infty$ are dense (Khudai-Veronov). 
\item If a set of leaves is measurable and has positive Lebesgue measure in $\P^2$, then it has full Lebesgue measure.
\item All singular points are hyperbolic.
\end{enumerate}

Note that if all the singularities are hyperbolic, then the leaves are uniformized by the unit disc, see Glutsyuk \cite{Glutsyuk} and Lins-Neto \cite{LinsNeto}.

Recall that foliations are given locally by holomorphic vector fields and their leaves are locally, integral curves of these vector fields. 
In the case of dimension 2 that we consider, foliations can be also given locally by a non-zero holomorphic 1-form and the leaves are Riemann surfaces on which this form vanishes. (Using rational vector fields, we see that projective complex surfaces admit large families of foliations. We will discuss later foliations on a compact K\"ahler surface).

If a holomorphic vector field $F$ has an isolated zero at some point $p$, we say that  the singularity $p$ is {\it hyperbolic} if the two eigenvalues of the linear part of the vector field at $p$ have non-real quotient. According to Poincar\'e, holomorphic vector fields are linearizable near hyperbolic singularities : if $p$ is such a singular point, then there are local holomorphic coordinates centered at $p$  such that the vector field has the form
$$z_2{\partial\over \partial z_2}-\eta z_1 {\partial\over\partial z_1}\ \raisepunct{,}$$
where $(z_1,z_2)\in\C^2$, $\eta=a+ib$, $a,b\in\R$, and $b\not=0$.

The statistical behavior of the leaves of a foliation can be studied using invariant positive $\ddc$-closed currents which are also called directed positive harmonic current.  Recall that $\dc:={i\over2\pi}(\dbar-\partial)$ and $\ddc = {i\over \pi}\ddbar$. Recall also
 that a positive current $T$ of bi-dimension $(1,1)$ is {\it directed by the foliation $\Fc$} if $T\wedge \Omega=0$ for every local holomorphic 1-form $\Omega$ defining $\Fc$. 
Let $\B$ be any flow box of $\Fc$ outside the singularities and denote by $V_\alpha$ the plaques of $\Fc$ in $\B$ parametrized by $\alpha$ in some transversal $\Sigma$ of $\B$. If $T$ is a positive $\ddc$-closed $(1,1)$-current directed by $\Fc$, then it has no mass on the singularities of $\Fc$ because this set is finite, see e.g. \cite{BS, Skoda}. On the flow box $\B$, this current 
 has the form
$$T_{|\B}=\int_{\alpha\in\Sigma} h_\alpha [V_\alpha] d\mu(\alpha),$$
where $h_\alpha$ is a positive harmonic function on $V_\alpha$. 

 In another language (harmonic measures), this was introduced by L. Garnett  for foliations without singularities in \cite{Garnett}, see also Sullivan \cite{Sullivan}. The existence of positive $\ddc$-closed currents directed by $\Fc$ was proved in the case with singularities by Berndtsson and the second author in \cite{BS}. When a leaf $L_w$ through a point $w$ is uniformized by the unit disc (i.e. it is hyperbolic in the sense of Kobayashi), an average on $L_w$ was introduced by  Forn\ae ss and the second author in \cite{FS1}. It allows us to get another construction of directed positive $\ddc$-closed currents. 
 
 More precisely, let $\D$ and $\D_r$ denote  the unit disc and the disc of center 0 and radius $r$ in $\C$. Let 
 $\phi^w:\D\to L_w$ be a universal covering map for $L_w$ with $\phi^w(0)=w$. 
Define the Nevanlinna characteristic  function for $\phi^w$ by 
$$T^w(r):=\int_0^r{dt\over t} \int_{\D_t} (\phi^w)^*(\omega_\FS), $$
 where  $\omega_\FS$ is the Fubini-Study form on $\P^2$. 
 Define the Nevanlinna current of index $r$, $0<r<1$, associated with $L_w$ by
 $$\tau_r^w :={1\over T^w(r)} (\phi^w)_*\Big[\log^+{r\over |\zeta |}\Big]={1\over T^w(r)} \int_0^r {dt\over t} (\phi^w)_*[\D_t].$$
Here, $\log^+:=\max(\log,0)$ and $\zeta$ is the standard coordinate of $\C$ so that the unit disc $\D$ is equal to $\{|\zeta|<1\}$. Note that for each $w$, the map $\phi^w$ is uniquely defined up to a rotation in $\D$. So the above definitions do not depend on the choice of $\phi^w$.

 It was shown in \cite{FS1} that for arbitrary isolated  singularities of $\Fc$ (not necessarily hyperbolic), $T^w(r)\to \infty$ as $r\to 1$. Therefore, the cluster points of $\tau^w_r$ are all $\ddc$-closed currents directed by $\Fc$. It turns out that a Birkhoff type theorem implies that for a generic foliation all extremal directed positive $\ddc$-closed currents of mass 1 can be obtained in this way \cite{DNS}, see also Nguyen \cite{Nguyen2}. General directed positive $\ddc$-closed currents are averages of the extremal ones.
 
 \medskip

Here is a main result of the present paper.

\begin{theorem} \label{th_main_P2}
Let $\Fc$ be a holomorphic foliation with only hyperbolic singularities in $\P^2$, seen as a compactification of $\C^2$. Assume that the only invariant algebraic curve is the line at infinity $L_\infty$. Then there is a unique positive $\ddc$-closed current of mass $1$ directed by $\Fc$ and it is the current of integration on $L_\infty$. In particular, if $\phi^w:\D\to L_w$ is a universal covering map of a leaf $L_w$ as above, then $\tau^w_r\to [L_\infty]$, in the sense of currents, as $r\to 1$. 
\end{theorem}

The result is surprising. Indeed, for generic foliations in $\Ac_d$, the leaves (except the line at infinity) are dense in $\P^2$. So one could expect that they spend a significant amount of hyperbolic time in every open set and that there should be a fat $\ddc$-closed non-closed current with support $\P^2$. Moreover, in the dynamics of holomorphic maps the canonical invariant currents and measures are not too singular; in particular, they give no mass to analytic sets, see the survey \cite{DS4}.

\medskip

Theorem \ref{th_main_P2} is in fact a particular case of the following result.

 \begin{theorem} \label{th_main_P2_bis}
 Let $\Fc$ be a holomorphic foliation by Riemann surfaces in $\P^2$ whose singularities are all hyperbolic. Then
 either (1) $\Fc$ has no invariant algebraic curve and then it has a unique positive $\ddc$-closed invariant $(1,1)$-current of mass 1; moreover, this unique current is not closed, or (2) $\Fc$ admits a finite number of invariant algebraic curves and any invariant positive $\ddc$-closed $(1,1)$-current is a combination of currents of integration on these curves; in particular, all invariant positive $\ddc$-closed $(1,1)$-currents are closed.
 \end{theorem}

The non-existence of a positive closed current on an exceptional minimal set, proved by Camacho-Lins Neto-Sad  in 
\cite{CLNS}, easily implies that for $\Fc$ as in the last theorem
the only positive {\bf closed} currents directed by $\Fc$ are currents of integration on algebraic leaves. 
When $\Fc$ has no invariant algebraic curve, we deduce that there is no non-zero positive closed $(1,1)$-current directed by $\Fc$. 
Forn\ae ss and the second author proved in \cite{FS2} that in this case, there is a unique positive $\ddc$-closed current of mass 1 directed by $\Fc$. So in order to obtain Theorem \ref{th_main_P2_bis}, we can assume that $\Fc$ admits an invariant algebraic curve $V$. The  key point of the proof is to show that if $T$ is a positive $\ddc$-closed current directed by $\Fc$ having no mass on any leaf, then $T$ is zero. 

For this purpose, we develop a theory of densities of positive $\ddc$-closed currents in a compact K\"ahler manifold with respect to a submanifold. The theory was developed by the authors in \cite{DS2} for positive closed currents and applications can be found in \cite{DNT, DS3}. Here, for simplicity, we limit ourselves to the case of positive $\ddc$-closed $(1,1)$-currents $T$ in a compact K\"ahler surface $X$.
Consider an irreducible smooth submanifold $V$ of dimension 1 of $X$. Assume for simplicity that $T$ has no mass on $V$.
We describe now the rough idea in this setting.   

In order to visualize the densities of $T$ along $V$, we dilate the coordinates in the normal direction to $V$. When the expansion factor tends to infinity, the image of $T$ by the dilation admits limit values which are positive closed $(1,1)$-currents $S$ on the normal line bundle $N_{V|X}$ of $V$ in $X$. We show that such currents $S$ are well-defined and we call them {\it tangent currents to $T$ along $V$}. We then show that $S$ has the form $\pi^*(\nu)$, where $\pi:N_{V|X}\to V$ is the canonical projection and $\nu$ is a positive measure on $V$. Moreover, the mass of $\nu$ is the cup-product $\{T\}\smallsmile \{V\}$ of the cohomology classes of $T$ and $V$. 

Now, in the case where $X=\P^2$, the cup-product  $\{T\}\smallsmile \{V\}$ vanishes only when $T=0$. Furthermore, if $T$ is directed by a foliation $\Fc$ as above and $V$ is an invariant curve, we can show that the tangent currents to $T$ along $V$ vanish. It follows immediately that $\{T\}\smallsmile \{V\}=0$ and then $T=0$ as desired. Both the theory of densities and the vanishing of the tangent currents to $T$ require delicate estimates. 

It is likely that if all singularities are linearizable, our analysis can be 
developed. However, it seems unlikely to get the uniqueness in general.
When some singularities are not linearizable, the local dynamics is already quite involved, and one should probably understand first harmonic currents near those singularities. This may be related to the Dulac-Moussu problem, see \cite{PM}.

The paper is organized as follows. In the next section, we give some basic properties of positive $\ddc$-closed currents and introduce the theory of densities for them. The proofs for the main results of the theory of densities are postponed to Section \ref{section_density_proof}. The last two sections contain some results for foliations on compact K\"ahler surfaces and, in particular, the proofs of Theorems \ref{th_main_P2} and \ref{th_main_P2_bis} stated above.

\medskip
\noindent
{\bf Acknowledgement.} The paper was partially written during the visit of the second  author 
at the National University of Singapore.  He would like to thank this university for its very warm hospitality and support.
We also thank the referee for his comments which helped to improve the exposition.

\section{Positive $\ddc$-closed currents on a compact K\"ahler surface} \label{section_currents}

The main tool used in this paper is the theory of densities for positive $\ddc$-closed currents. As mentioned above, for simplicity, we  only consider the case of $(1,1)$-currents in a compact K\"ahler surface, and for the reader's convenience, a self-contained exposition will be given. 
Let $X$ be a compact K\"ahler surface. We start this section with some basic properties of positive $\ddc$-closed currents, see e.g. \cite{BS, Demailly, DS1, Skoda} and also  \cite{Voisin} for compact K\"ahler manifolds.

Let $T$ be a $\ddc$-closed $(1,1)$-current on $X$. Let $\alpha$ and $\alpha'$ be two smooth closed $(1,1)$-forms on $X$ which are in the same cohomology class in $H^{1,1}(X,\C)$. The $\ddc$-lemma implies the existence of a smooth function $f$ such that $\alpha-\alpha'=\ddc f$. It follows from Stokes formula that 
$$\langle T,\alpha\rangle -\langle T,\alpha'\rangle =\langle T,\ddc f\rangle =\langle \ddc T,f\rangle =0.$$
So $T$ defines an element in the dual space of $H^{1,1}(X,\C)$. By Poincar\'e's duality theorem, the dual of $H^{1,1}(X,\C)$ is naturally identified with itself. Therefore, $T$ defines a class in $H^{1,1}(X,\C)$ that we denote by $\{T\}$. 
When $T$ is moreover a real current, e.g. when $T$ is positive, this class belongs to $H^{1,1}(X,\R)$ which, by definition, is the intersection $H^{1,1}(X,\C)\cap H^2(X,\R)$. 
We have the following basic result.

\begin{proposition} \label{prop_T_V}
Let $V$ be an irreducible analytic subset of dimension $1$ of $X$ and
let $T$ be a positive $\ddc$-closed $(1,1)$-current on $X$.  Then there is a unique constant $c\geq 0$ such that $T=c[V]+S$, where $[V]$ is the current of integration on $V$ and $S$ is a positive $\ddc$-closed current without mass on $V$. 
\end{proposition}
\proof
Let $S$ denote the restriction of $T$ to $X\setminus V$. This is a positive $(1,1)$-current on $X\setminus V$ with finite mass. Therefore, it can be extended by 0 to a current on $X$ that we still denote by $S$. We have $\ddc S\leq 0$, see e.g. \cite{AB, DS1}. However, by Stokes theorem, we have 
$$\int_X\ddc S= \langle \ddc S, 1\rangle = \langle S,\ddc 1\rangle =\langle S,0\rangle =0.$$
Hence, $\ddc S=0$. The $(1,1)$-current $T-S$ is positive $\ddc$-closed and supported by $V$. Hence, it is equal to $c[V]$, where $c$ is a positive harmonic function on $V$ and $[V]$ is the current of integration on $V$. Since $V$ is compact,  $c$ is necessarily constant. The proposition follows.
\endproof

We will need later the following lemma.

\begin{lemma} \label{lemma_harm_c11}
Let $K$ be a compact subset of $X$ and  $T$ be a positive $\ddc$-closed $(1,1)$-current on $X$ without mass on $K$. Let $\Psi:=\alpha+\ddc u$, where $\alpha$ is a smooth closed $(1,1)$-form on $X$ and $u$ is a function on $X$ which is smooth outside $K$ and of class $\Cc^{1,1}$ on $X$, i.e. its partial derivatives of order $1$ are Lipschitz functions. Then the integral $\langle T,\Psi\rangle_{X\setminus K}$ of $T\wedge \Psi$ on $X\setminus K$  is equal to $\{T\}\smallsmile \{\alpha\}$. 
\end{lemma}
\proof
Consider first the case where $u$ is smooth. Since $T$ has no mass on $K$, we have  $\langle T,\Psi\rangle_{X\setminus K}=\langle T,\Psi\rangle$. Since $\Psi$ is cohomologous to $\alpha$, the above discussion on the cohomology classes of $\ddc$-closed currents implies that the last integral is equal to $\{T\}\smallsmile \{\alpha\}$ and the lemma follows.

Now, we will reduce the general case to the first case using a partition of unity and convolutions on the charts of $X$. There is a sequence of smooth functions $u_n$ on $X$ such that 
\begin{itemize}
\item[(1)] $u_n$ and its partial derivatives of order $\leq 2$ are bounded uniformly on $n$;
\item[(2)] $u_n$ converges to $u$, together with all partial derivatives of order $\leq 2$, locally uniformly on $X\setminus K$. 
\end{itemize}
Therefore, if $\Psi_n:=\alpha+\ddc u_n$, we have since $T$ has no mass on $K$
$$\langle T,\Psi\rangle_{X\setminus K}=\lim_{n\to \infty} \langle T,\Psi_n\rangle = \lim_{n\to\infty}\{T\}\smallsmile \{\alpha\}=\{T\}\smallsmile \{\alpha\}.$$
This completes the proof of the lemma.
\endproof

We continue our discussion by giving a construction of positive $\ddc$-closed currents. 
Let $\phi:\D\to X$ be a holomorphic map from the unit disc $\D$ to $X$. 
We say that $\phi$ is a parametrized holomorphic disc. Fix a K\"ahler form $\omega$ on $X$.
For $0<r<1$, define {\it the Nevanlinna characteristic function} $T^\phi(r)$ by 
$$T^\phi(r):=\int_0^r{dt\over t} \int_{\D_t} \phi^*(\omega).$$
The mass of the averaging current
$$\tau_r^\phi :={1\over T^\phi(r)} \phi_*\Big[\log^+{r\over |\zeta |}\Big]={1\over T^\phi(r)} \int_0^r {dt\over t} \phi_*[\D_t]$$
is 1. Note that this current is positive and a direct computation gives
$$\ddc \tau_r^\phi={1\over T^\phi(r)}\phi_*(\nu_r-\delta_0),$$
where $\nu_r$ is the Haar probability measure on the boundary of $\D_r$ and $\delta_0$ is the Dirac mass at $0\in\D$.

\begin{definition} \rm
We call {\it Nevanlinna current associated with $\phi$} any cluster value of $\tau_r^\phi$ when $r\to 1$. 
The integral 
$$T^\phi(1):=\lim_{r\to 1} T^\phi(r)=\int_0^1{dt\over t} \int_{\D_t} \phi^*(\omega)$$
is called {\it the Nevanlinna area} of the parametrized holomorphic disc $\phi$.
\end{definition}

Observe that the above computation on $\ddc\tau_r^\phi$ shows that  when $T^\phi(1)$ is infinite, the Nevanlinna currents are  $\ddc$-closed. 

We are interested in positive $\ddc$-closed currents directed by a holomorphic foliation on $X$. For a generic foliation, by a recent ergodic theorem, such currents can be obtained as averages of Nevanlinna currents associated with a natural parametrization of the leaves of the foliation, see \cite{DNS}. 

\begin{definition} \rm
 Let $U$ be an open subset of $X$. We say that the map {\it $\phi$ crosses $U$ properly}  if the restriction of $\phi$ to each connected component of $\phi^{-1}(U)$ is a proper map into $U$. We say that {\it $\phi$ crosses $U$ locally properly} if every point in $U$ admits a neighbourhood $U'\subset U$ such that $\phi$ crosses $U'$ properly.
 \end{definition}

Note that if $\phi$ crosses $U$ properly, it crosses properly  all open subsets of $U$.
It is clear that if $\phi$ is a universal covering of a leaf of a foliation, then it crosses any regular flow box properly.  
Recall that a set  $E$ is {\it complete pluripolar} if  every $p\in X$ admits a neighbourhood $U_p$ such that $E\cap U_p=\{u=-\infty\}$ for some p.s.h. function $u$ on $U_p$. We have the following criterion for the non-finiteness of the Nevanlinna area.

\begin{proposition} \label{prop_current_Nev}
Let $E\subset X$ be a closed complete pluripolar subset of $X$ and $\phi:\D\to X\setminus E$ be a parametrized holomorphic disc. Assume that $\phi$ crosses $X\setminus E$ locally properly. Then the Nevanlinna area of $\phi$ is infinite.
\end{proposition}
\proof
Observe that since $\phi$ crosses $X\setminus E$ locally properly, it is not a constant map. 
Assume that the Nevanlinna area of $\phi$ is finite. Then,  $T:=\phi_*(-\log|\zeta|)$ is a positive $(1,1)$-current of finite mass in $X$. The hypothesis implies that it is $\ddc$-closed outside $E\cup\{\phi(0)\}$. Hence, $\ddc T$ is negative on $X$, see e.g. \cite[Th. 1.3]{DS1}. Since $X$ is compact,  using Stokes theorem as in Proposition \ref{prop_T_V}, we get  $\ddc T=0$. 
On the other hand, there is a neighbourhood $U$ of $\phi(0)$ such that $\phi$ crosses $U$ properly. We then deduce that $\ddc T$ contains the Dirac mass at $\phi(0)$. 
This is a contradiction.  
\endproof

\begin{proposition} \label{prop_current_leaf}
Let $E$ and $\phi$ be as in the last proposition. Assume that there exists a non-zero positive $\ddc$-closed $(1,1)$-current $T$ on $X$ without mass on $E$ and on  $X\setminus \phi(\D)$. Then $\overline {\phi(\D)}$ is an irreducible analytic subset of dimension $1$ of $X$ and $T$ is equal to a constant times the current of integration on it.  
\end{proposition}
\proof
Consider a point in $X\setminus E$. There is a neighbourhood $U$ of  this point such that $\phi$ crosses $U$ properly. It follows that the intersection $\phi(\D)\cap U$ is a finite or countable union of irreducible analytic subsets of $U$. Denote them by $V_1,V_2,\ldots$. The restriction of $T$ to each plaque $V_i$ has the form $h_i [V_i]$, where $h_i$ is a positive subharmonic function on $V_i$, see e.g. \cite[Th. 1.3]{DS1}. It follows from the hypotheses that the restriction of $T$ to $U$ is equal to $\sum_i h_i[V_i]$.  Since $T$ is harmonic, we deduce that $h_i$ is harmonic for all $i$. 

Observe that $\phi(\D)$ is an irreducible Riemann surface immersed in $X\setminus E$. We deduce that $T$ has the form $h[\phi(\D)]$, where $h$ is a positive harmonic function on $\phi(\D)$, i.e. $\Delta h=0$ in each plaque of $\phi(\D)$. Consider for any constant $c\geq 0$ the current $T_c:=\min(c,h)[\phi(\D)]$. This is a positive $(1,1)$-current on $X$ having no mass on $E$ and  $\ddc T_c\leq 0$ on $X\setminus E$. Therefore, we have $\ddc T_c\leq 0$ on $X$, see \cite{AB, DS1}. Using Stokes theorem, we obtain that $\ddc T_c=0$ on $X$ and hence $\min(c,h)$ is harmonic on $\phi(\D)$. The property holds for any $c\geq 0$. It is not difficult to deduce that $h$ is constant. It is strictly positive because $T\not=0$.

So $[\phi(\D)]$ is a positive current of finite mass on $X$ or equivalently the area of $\phi(\D)$ is finite. We deduce that the intersection of $\phi(\D)$ with each open set $U$ as above is a locally finite  union of analytic sets of $U$. Thus, $\phi(\D)$ is an irreducible analytic subset of $X\setminus E$ of finite area. By a Bishop type theorem, $\overline{\phi(\D)}$ is an analytic subset of $X$ which is irreducible because  $\phi(\D)$ is irreducible, see also \cite{DS1, Skoda}. We also obtain that $T$ is equal to a constant times the current of integration on $\overline{\phi(\D)}$.
\endproof

From now on, we assume that  $V$ is a smooth irreducible analytic subset of dimension 1 of $X$.  So $V$ is a smooth connected compact Riemann surface.  Consider a positive $\ddc$-closed current $T$ on $X$ without mass on $V$. We will study the density of $T$ near $V$ via a notion of "tangent cone" to $T$ along $V$ that we introduce now. 

Denote by $\Tan(X)$  the tangent bundle of $X$ and $\Tan(X)_{|V}$ its restriction to $V$. Let $\Tan(V)$ be the tangent bundle of $V$ which is naturally identified with a sub-bundle of $\Tan(X)_{|V}$. The quotient $\Tan(X)_{|V}/\Tan(V)$ is the normal bundle of $V$ in $X$ and will be denoted by $N_{V|X}$. We identify $V$ with the zero section in each vector bundle $\Tan(X)_{|V}$, $\Tan(V)$ and $N_{V|X}$.
 
\begin{definition} \rm
A {\it smooth admissible map} is a smooth  bijective map $\tau$ from a neighbourhood of  $V$ in $X$  to a neighbourhood of $V$ in $N_{V|X}$ such that 
\begin{enumerate}
\item The restriction of $\tau$ to $V$ is the identity map from $V$ to $V$; in particular, the restriction of the differential $d\tau$ to $V$ induces three endomorphisms of $\Tan(X)_{|V}$, $\Tan(V)$ and $N_{V|X}$ respectively; 
\item The differential $d\tau(x)$, at any point $x\in V$, is a holomorphic endomorphism of the tangent space to $X$ at $x$;
\item The endomorphism of $N_{V|X}$, induced by  $d\tau$,  is  the identity map.
\end{enumerate}
\end{definition}
Note that the dependence of $d\tau(x)$ in $x\in V$ is not necessarily holomorphic. When $X=\P^2$ and $V$ is a projective line, one can construct such a map which is holomorphic. Indeed, we consider $\P^2$ as a  compactification of $\C^2$ with $V$ as the line at infinity. Then the normal bundle $N_{V|X}$ can be identified with $\P^2\setminus\{0\}$ and its fibers are identified with the lines through 0. With this identification, the identity map is holomorphic and admissible. If $V$ is a general curve in $\P^2$, there is no admissible holomorphic map. In the general case, a smooth admissible map can be constructed in the following way. 

\begin{example} \rm \label{ex_admissible}
Fix a smooth K\"ahler metric on $X$ associated with a K\"ahler form $\omega$ (here, it is enough to consider a smooth Hermitian metric but we will need later in Proposition \ref{prop_ex_admissible} that it is K\"ahler). For each point $x\in V$ denote by $N_x$ the orthogonal complement of the tangent line to $V$ at $x$ in the tangent plane to $X$ at $x$, with respect to the considered metric. The union of $N_x$ for $x\in V$ can be identified with $N_{V|X}$ but this identification is not a holomorphic map in general. We construct the map $\tau^{-1}$ from a neighbourhood of $V$ in $N_{V|X}$ to a neighbourhood of $V$ in $X$ in the following way: for $z\in N_x$ close enough to $x$, $\tau^{-1}(z)$ is the image of $z$ by the exponential map, which is defined on the tangent plane to $X$ at $x$. It is easy to check that $\tau$ is well-defined and is smooth admissible with $d\tau(x)=\id$ for $x\in V$.  We don't define general admissible maps as maps with $d\tau =\id$ on $V$ because we will use a larger class.
\end{example}

For $\lambda\in \C$, denote by $A_\lambda$ the map from the line bundle $N_{V|X}$ to itself induced by the multiplication by $\lambda$ in each fiber. Let $\tau$ be a smooth admissible map as above. Define 
$$T_\lambda:=(A_\lambda)_*\tau_*(T).$$
This is a current of degree 2. Its domain of definition is  some open subset of $N_{V|X}$ containing $V$ which increases to $N_{V|X}$ when $|\lambda|$ increases to infinity. Note that $T_\lambda$ is not a $(1,1)$-current and we cannot speak of its positivity. Moreover, it is not $\ddc$-closed in general and we cannot speak of its cohomology class. When $T$ is a positive closed current, then $T_\lambda$ is closed and the situation is simpler. 
Denote by $\pi: N_{V|X}\to V$ the canonical projection. We have the following theorem.

\begin{theorem} \label{th_tangent}
The mass of $T_\lambda$ on any given compact subset of $N_{V|X}$ is bounded when $|\lambda|$ is large enough. If $S$ is a cluster value of $T_\lambda$ when $\lambda\to\infty$, then it is a positive closed $(1,1)$-current on $N_{V|X}$  given by 
$S=\pi^*(\nu),$
where $\nu$ is a positive measure on $V$ of mass equal to $\{T\}\smallsmile \{V\}$. Moreover, if $(\lambda_n)$ is a sequence tending to infinity such that $T_{\lambda_n}\to S$, then $S$ may depend on $(\lambda_n)$ but it does not depend on the choice of the map $\tau$. 
\end{theorem}

Note that if the class $\{V\}$ of $V$ is K\"ahler (e.g. when $X=\P^2$), we have $\{T\}\smallsmile \{V\}=0$ only when $T=0$. So in this case, $S=0$ implies $T=0$. We will use this property later. 
The proof of the above theorem will be given in the next section. We can now introduce the following notion.

\begin{definition} \rm
Any current $S$ obtained as in Theorem \ref{th_tangent} is called {\it a tangent current} to $T$ along the curve $V$.
\end{definition}

Note that in general $S$ is not unique as this is already the case for positive closed currents, see \cite{DS2} for details. 

In the rest of this section, we will explain how to compute tangent currents using local coordinates. Fix a local holomorphic coordinate system $z=(z_1,z_2)$ of $X$, with $|z_1|<3$ and $|z_2|<3$,  defined on an open set  such that $V$ is given there by the equation $z_2=0$. Consider the open set $\U:=\{|z_1|<1, |z_2|<1\}$ of $X$ and for simplicity, we identify it with the unit bi-disc $\D^2$ in $\C^2$. 
So $V\cap \D^2$ is identified with $\D\times \{0\}$.  

We can also, in a natural way, identify $N_{V|X}$ with $\D\times\C$ which is an open subset of $\C^2$ containing $\D^2$. 
In these local coordinates, the map $A_\lambda$, introduced above, is given by 
$$a_\lambda(z_1,z_2)=(z_1,\lambda z_2).$$
We use here the notation $a_\lambda$ in order to avoid the confusion with the global map $A_\lambda$ on $N_{V|X}$. 
Tangent currents can be computed locally using the following result.

\begin{proposition} \label{prop_tangent_local}
Let $S$ and $\lambda_n$ be as in Theorem \ref{th_tangent}. Then in the above local coordinates we have 
$$S=\lim_{n\to\infty} (a_{\lambda_n})_*(T) \quad \text{on}\quad \D\times \C.$$
\end{proposition}

We also postpone the proof of this proposition to the next section. Note that this proposition also shows that $S$ does not depend on the choice of $\tau$ because the limit in that proposition does not involve $\tau$.

\section{Existence and properties of tangent currents}  \label{section_density_proof}

In this section, we will prove the last two results from the previous section. 
We use the notations introduced before Proposition \ref{prop_tangent_local}. The following two lemmas are the main technical results. 

\begin{lemma} \label{lemma_technique_0}
For each $0<r\leq 1$, there are smooth positive closed $(1,1)$-forms $S_n$ on $X$, $n\geq 0$, 
such that 
$$\sum_{n\geq 0} \|S_n\| \leq c \qquad \text{and} \qquad  ir^{-2} dz_2\wedge d\overline z_2\leq \sum_{n\geq 0} S_n\quad \text{on}\quad \{|z_1|\leq  1, 0<|z_2|\leq r\},$$
where $c>0$ is a constant independent of $r$. 
\end{lemma}
\proof
Let $m\geq 0$ be the integer such that $e^{-m-1}<r\leq e^{-m}$. We have $r^{-2}\leq e^{2m+2}$. 
We will construct  smooth closed $(1,1)$-forms $R_m$ of mass bounded by a constant such that $R_m\geq ie^{2m}dz_2\wedge d\overline z_2$ on $\{|z_1|\leq  1, e^{-m-1}\leq |z_2|\leq  e^{-m}\}$. The lemma follows immediately by taking $S_n:=e^{2-2n} R_{m+n}$ for $n\geq 0$.  

Let $\phi$ be a quasi-psh function on $X$ such that  $\ddc\phi=[V]-\alpha$ for some real smooth closed $(1,1)$-form $\alpha$. In the local coordinates $(z_1,z_2)$, we have $[V]=\ddc \log|z_2|$. It follows that $\psi:=\phi-\log|z_2|$ is a local potential of $-\alpha$ and therefore it is a smooth function on the domain $\{|z_1|<3, |z_2|<3\}$. 
Let $M\geq 1$ be a constant such that $|\psi|\leq M$ on $\{|z_1|\leq 2, |z_2|\leq 2\}$. 

Let $\chi:\R\to \R$ be an increasing convex smooth function such that $\chi(t)=0$ for $t\leq -3M$, $\chi(t)=t$ for $t\geq 3M$, $0\leq\chi'(t)\leq 1$,  and $\chi''(t)\geq 1/(10M)$  in $(-2M,2M)$. Fix also a K\"ahler form $\omega$ on $X$ and choose a constant $A\gg M$ large enough. Define
$$R_m:=A\ddc \big[\chi(\phi+m)\big] +A^2\omega.$$ 
This is clearly a smooth closed $(1,1)$-form on $X$. We first show that it is positive and has bounded mass.

A direct computation gives
$$R_m=A\chi'(\phi+m)\ddc \phi + {A\over \pi} \chi''(\phi+m)i \partial\phi\wedge \dbar\phi + A^2\omega.$$
The second term is positive. The first term is bounded below by $-Ac\omega$ if $\alpha\geq -c\omega$. We then deduce that $R_m$ is positive since $A$ is chosen large enough. Now, since $R_m$ is cohomologous to $A^2\omega$, its mass is equal to the mass of $A^2\omega$ and hence is bounded.

It remains to check that $R_m\geq e^{2m} idz_2\wedge d\overline z_2$ on the domain $\{|z_1|\leq 1, e^{-m-1}\leq |z_2| \leq e^{-m}\}$. 
We need to consider the second term in the above expansion of $R_m$. 
In the considered domain, we have $|\phi+m|\leq 2M$ and therefore $\chi''(\phi+m)\geq 1/(10M)$. 
On the other hand, the form  $i \partial\phi\wedge \dbar\phi$ is equal to
\begin{eqnarray*}
\lefteqn{i \partial(\psi+\log|z_2|)\wedge \dbar(\psi+\log|z_2|)} \\
& = & -3i \partial \psi\wedge \dbar \psi +  i \partial\Big(2\psi+{1\over 2}\log|z_2|\Big)\wedge \dbar\Big(2\psi+{1\over 2}\log|z_2|\Big) +{3\over 4} i\partial \log|z_2|\wedge \dbar \log |z_2| \\
& \geq & -3i \partial \psi\wedge \dbar \psi + {3\over 16} e^{2m} idz_2\wedge d\overline z_2 \quad \text{since the second term in the last line is positive.}
\end{eqnarray*}
The first term in the last line is bounded. Now, we can easily deduce from the above expansion of $R_m$ the desired inequality for $A\gg M$ large enough.
\endproof

Recall that $T$ is a positive $\ddc$-closed $(1,1)$-current without mass on $V$. 

\begin{lemma} \label{lemma_technique}
There is a constant $c>0$ and for each $r$, $0<r<1$, a constant  $\epsilon_r>0$ depending on $T$, 
such that $\epsilon_r\to 0$ as $r\to 0$ and such that the following estimates hold. For any continuous function $f(z)$ with compact support in $\{|z_1|< 1, |z_2|< r\}$ with $|f|\leq 1$, we have
$$|\langle T, f(z) dz_1\wedge d\overline z_1\rangle| \leq \epsilon_r, \qquad |\langle T, f(z) dz_2\wedge d\overline z_2\rangle| \leq c r^2$$
and
$$|\langle T, f(z) dz_1\wedge d\overline z_2\rangle |\leq \epsilon_r r,
 \qquad  |\langle T, f(z) d\overline z_1\wedge d z_2\rangle |\leq \epsilon_r r.$$
\end{lemma}
\proof
Since $T$ has no mass on $V$, the measure $T\wedge idz_1\wedge d\overline z_1$ has no mass on $V$. It follows that its mass on $\{|z_1|\leq 1, |z_2|\leq r\}$ tends to 0 as $r\to 0$. The first estimate follows.

For the second estimate, it is enough to consider the case where $f$ is positive. So $r^{-2}f(z) idz_2\wedge d\overline z_2$ is a positive form. By Lemma \ref{lemma_technique_0},  since $T$ has no mass on $V$, we have 
$$ |\langle T, r^{-2} f(z) dz_2\wedge d\overline z_2\rangle | \leq  \sum_{n\geq 0}\langle T, S_n \rangle =\sum_{n\geq 0} \{T\}\smallsmile \{S_n\}.$$
The last expression is bounded because the sum of the masses of $S_n$ is bounded and hence the sum of their  cohomology classes is also bounded. This proves the second estimate. Note that it is crucial here that $T$ is a global positive $\ddc$-closed current and the $S_n$'s are global positive closed forms. 

We now prove the third inequality.  Let $\chi$ be a smooth function with compact support in 
$\{|z_1|< 1, |z_2|< r\}$ such that $0\leq\chi\leq 1$ and $\chi=1$ in a neighbourhood of the support of $f$.
We have by the Cauchy-Schwarz inequality
\begin{eqnarray*}
|\langle T, f(z) dz_1\wedge d\overline z_2\rangle | & = & |\langle T,\chi(z) f(z) dz_1\wedge d\overline z_2\rangle |\\
& \leq &  |\langle T, \chi(z)^2 dz_1\wedge d\overline z_1\rangle |^{1/2} |\langle T, f(z)^2 dz_2\wedge d\overline z_2\rangle |^{1/2}.
\end{eqnarray*}
We deduce the third inequality in the proposition from the first  and second ones by replacing $\epsilon_r$ with $\sqrt{c\epsilon_r}$. 
The last inequality can be obtained in the same way.
\endproof

\begin{definition} \rm \label{def_negligible}
Let $(\alpha_\lambda)$ be a family of  2-forms in $X$ or in $N_{V|X}$ depending on $\lambda\in \C$ with $|\lambda|$ larger than a constant. We say that this family is {\it negligible} if the support $\supp(\alpha_\lambda)$ tends to $V$ as $\lambda\to\infty$ and if we have
in any local coordinate system as above
\begin{itemize}
\item[(1)] $\supp (\alpha_\lambda)\cap \{|z_1|<1, |z_2|<1\}$ is contained in $\{|z_2|\leq A|\lambda|^{-1}\}$ for some constant $A>0$ independent of $\lambda$;  
\item[(2)] The coefficient of $dz_1\wedge d\overline z_1$ in $\alpha_\lambda$ is bounded, the sup-norm of the coefficient of $dz_2\wedge d\overline z_2$ is $o(\lambda^2)$ and 
the sup-norms of the other coefficients are $O(\lambda)$ as $\lambda\to\infty$.  
\end{itemize}
\end{definition}

Note that Property (1) is often easy to check. If we use the coordinates $(z_1,\lambda z_2)$ instead of $(z_1,z_2)$, then Property (2) is equivalent to (2a) the coefficients of $\alpha_\lambda$ are bounded and (2b) the sup-norm of the coefficient of $i d(\lambda z_2)\wedge d(\overline\lambda\overline z_2)$ tends to 0 as $\lambda\to\infty$. Property (2a) is often easy to check. So to check that a family is negligible, one often only has to bound the coefficient of $idz_2\wedge d\overline z_2$. 

We will need the notion of negligible family in the study of tangent currents. More precisely, the following lemma will be used later.

\begin{lemma} \label{lemma_negligible}
Let $(\alpha_\lambda)$ be a negligible family of smooth $2$-forms in $X$. Let $T$ be a positive $\ddc$-closed current on $X$ without mass on $V$. 
Then $\langle T,\alpha_\lambda\rangle \to 0$ as $\lambda\to\infty$. 
\end{lemma}
\proof
We can use a partition of unity in order to work with local coordinates. So we can assume that the forms $\alpha_\lambda$ 
have supports in $\{|z_1|<1,|z_2|<1\}$. Since $T$ is of bidegree $(1,1)$, we only need to test forms of bidegree $(1,1)$. 
Lemma \ref{lemma_technique}, applied to $r:=A|\lambda|^{-1}$ with $A$ from Definition \ref{def_negligible}, shows that   $\langle T,\alpha_\lambda\rangle \to 0$.
\endproof

We need a description of $\tau$ in local coordinates $(z_1,z_2)$. Consider the Taylor  expansion of order 2 of $\tau$ in $z_2,\overline z_2$ with functions in $z_1$ as coefficients. 
Since $\tau$ is smooth admissible, when $z_2\to 0$, 
this map and its differential can be written as
$$\tau(z_1,z_2)= \Big(z_1+a(z_1)z_2 + O(|z_2|^2), z_2+ O(|z_2|^2)\Big)$$
and
 $$d\tau(z_1,z_2)= \Big(dz_1+O(1)dz_2+O(|z_2|)+O^*(|z_2|^2), dz_2+O^*(|z_2|^2)\Big),$$
where $a(z_1)$ is a smooth function in $z_1$ and $O^*(|z_2|^k)$ is any smooth 1-form that can be written as
$$O^*(|z_2|^k)=O(|z_2|^{k-1})dz_2+O(|z_2|^{k-1} )d\overline z_2+O(|z_2|^k).$$ 
Here, $O(|z_2|^k)$ denotes any smooth function or form whose sup-norm is bounded by a constant times $|z_2|^k$. 
We also have
 $$d\tau^{-1}(z_1,z_2)= \Big(dz_1+O(1)dz_2+O(|z_2|)+O^*(|z_2|^2), dz_2+O^*(|z_2|^2)\Big).$$

\begin{lemma} \label{lemma_tau_negli}
If $(\alpha_\lambda)$ is a negligible family of $2$-forms on $N_{V|X}$, then $(\tau^*(\alpha_\lambda))$ is also a negligible family of $2$-forms on $X$.
\end{lemma}
\proof
This is a direct consequence of the above local description of $d\tau$. 
\endproof

\begin{proposition} \label{prop_mass}
Let $\Phi$ be a continuous $2$-form with support in a fixed  compact subset $K$ of $N_{V|X}$. 
Define $\Phi_\lambda:=A_\lambda^*(\Phi)$ and $\Psi_\lambda:=\tau^*A_\lambda^*(\Phi)$. 
The following properties hold.
\begin{itemize}
\item[(1)]   If $\Phi\wedge \pi^*(\Omega)=0$ for any smooth $(1,1)$-form $\Omega$ on $V$, then the families of $(\Phi_\lambda)$ and $(\Psi_\lambda)$ are  negligible.
\item[(2)]  If $\|\Phi\|_\infty\leq 1$, then $\limsup_{\lambda\to\infty}|\langle T_\lambda,\Phi\rangle|$ is bounded by a constant which does not depend on $\Phi$.
\item[(3)]  If $\Phi\wedge \pi^*(\Omega)\geq 0$ for any smooth positive $(1,1)$-form $\Omega$ on $V$, then any limit value of 
$\langle T_\lambda,\Phi\rangle$, when $\lambda\to\infty$, is non-negative. In particular, the property holds when $\Phi$ is a positive $(1,1)$-form.
\item[(4)] If $\Phi=\ddc u$ for some smooth function $u$ with compact support in $N_{V|X}$, then $\langle T_\lambda,\Phi\rangle \to 0$ as $\lambda\to\infty$.
\end{itemize}
\end{proposition}
\proof
(1) If $(\chi_k)$ is a finite partition of unity for $V$, then $(\chi_k\circ \pi)$ is a finite partition of unity for $N_{V|X}$. Using such a partition,  we can reduce the problem to the case where  $\Phi$ has support in $\{|z_1|<1/2, |z_2|<A/2\}$ with $(z_1,z_2)$ as above and $A>0$ a constant. The hypothesis in (1) implies that 
the coefficient of $dz_2\wedge d\overline z_2$ in $\Phi$ vanishes. Then, a direct computation shows that $(\Phi_\lambda)$ is negligible. By Lemma \ref{lemma_tau_negli}, the family $(\Psi_\lambda)$ is also negligible.

\medskip

(2) As above, we can assume that  $\Phi$ has support in $\{|z_1|<1/2, |z_2|<A/2\}$. 
Modulo a negligible family of forms, thanks to the first assertion, we have $\Phi_\lambda\simeq f_\lambda(z) |\lambda|^2 idz_2\wedge d\overline z_2$, where $f_\lambda$ is a smooth function supported by $\{|z_1|<1/2, |z_2|<A|\lambda|^{-1}/2\}$ and $|f_\lambda|$ is bounded by a constant. Then, we deduce from the above expansion of $d\tau$ that $\Psi_\lambda$ satisfies a similar property with support in $\{|z_1|<1, |z_2|<A|\lambda|^{-1}\}$ when $\lambda$ is large enough. 
By Lemma \ref{lemma_negligible}, negligible families of forms do not change the limit we are considering. The second estimate in Lemma \ref{lemma_technique} implies the result.

\medskip

(3) We can assume that $\Phi$ is as in (2). The hypothesis of (3) implies that the coefficient of $idz_2\wedge d\overline z_2$ in $\Phi$ is positive. It follows that $f_\lambda\geq 0$. We also see using the expansion of $\tau$ that $\Psi_\lambda$ is the product of a positive function $g_\lambda$ with $idz_2\wedge d\overline z_2$ plus a form in a negligible family. Since $T$ is positive, we have $\langle T,g_\lambda idz_2\wedge d\overline z_2\rangle\geq 0$.  The result follows easily.

\medskip

(4) Using a partition of unity, we can assume that $u$ has support in $\{|z_1|<1/2, |z_2|<A/2\}$. Define $u_\lambda:=A_\lambda^*u=u\circ A_\lambda$ and 
$v_\lambda:=\tau^*A_\lambda^*u=u_\lambda\circ\tau$. We have 
$u_\lambda(z_1,z_2)=u(z_1,\lambda z_2)$ and $v_\lambda=u(\tau_1,\lambda\tau_2)$ if we write $\tau=(\tau_1,\tau_2)$ in the considered coordinates. Using the above local description of $\tau$ and $d\tau$, we easily check that 
$$\tau^*(\ddc u_\lambda)-\ddc v_\lambda$$
is in a negligible family of 2-forms. For example, when we expand $\ddc v_\lambda=\ddc (u_\lambda\circ \tau)$, we obtain the term 
$$\lambda {\partial u\over \partial z_2}(\tau_1,\lambda\tau_2)\ddc \tau_2=-\lambda {\partial u\over \partial z_2}(\tau_1,\lambda\tau_2)\dc (d\tau_2).$$ 
The derivative $\partial u/\partial z_2$ being bounded, the identity $d\tau_2 =dz_2+O^*(|z_2|^2)$ implies that this term is negligible.

Now, it follows from Lemma \ref{lemma_negligible} that 
$$\langle  T_\lambda,\ddc u\rangle =\langle T,\tau^*(\ddc u_\lambda)\rangle = \langle T,\ddc v_\lambda\rangle +o(1)$$
as $\lambda\to \infty$. On the other hand, since $T$ is $\ddc$-closed, we have $\langle T,\ddc v_\lambda\rangle =0$. The result follows. 
\endproof

Consider any sequence $(\lambda'_n)$ in $\C$ tending to infinity. 
The second assertion in Proposition \ref{prop_mass} implies that  we can extract a subsequence $(\lambda_n)$ such that $T_{\lambda_n}$ converges to a 2-current $S$ of locally finite mass in $N_{V|X}$.  The first assertion in that proposition shows that in the above local coordinates, if the coefficient of $dz_2\wedge d\overline z_2$ in $\Phi$ vanishes then $\langle S,\Phi\rangle=0$. 
It follows that  $S\wedge dz_1=0$ and $S\wedge d\overline z_1=0$, or equivalently,
 $S\wedge \pi^*(\Omega)=0$ for any smooth form $\Omega$ of degree $1$ or $2$ on $V$. 
Hence $S$ is a current of bi-degree $(1,1)$.  The third assertion of the last proposition implies that $S$ is positive. Finally, the fourth assertion is equivalent to saying that $S$ is $\ddc$-closed.

\begin{lemma} \label{lemma_decomp}
There is a positive measure $\nu$ on $V$ such that $S=\pi^*(\nu)$. In particular, $S$ is closed.
\end{lemma}
\proof
Consider the family $\Fc$ of all positive $\ddc$-closed $(1,1)$-currents $R$ on $N_{V|X}$ which are vertical in the sense that 
$R\wedge \pi^*(\Omega)=0$ for any smooth form $\Omega$ of degree $1$ or $2$ on $V$.

\medskip
\noindent
{\bf Claim. }
If $S$ is any current in $\Fc$ and $u$ is a smooth positive function on $V$, then $(u\circ\pi)S$ also belongs to $\Fc$. 

\medskip

Indeed, it is clear that $(u\circ\pi)S$ is positive and vertical. 
The only point to check is that 
 $(u\circ\pi)S$ is $\ddc$-closed. Define $\widetilde u:=u\circ\pi$. We have $\ddc S=0$ and since $S$ is vertical, we have $d\widetilde u\wedge S=0$, $d^c\widetilde u\wedge S=0$ and $\ddc \widetilde u\wedge S=0$. Therefore, a direct computation gives
 $$\ddc (\widetilde u S)= d(\dc\widetilde u\wedge S)-\dc(d\widetilde u\wedge S) -\ddc\widetilde u\wedge S +\widetilde u\ddc S=0,$$
 which completes the proof of the claim.

\smallskip

It follows from the claim that every extremal element in $\Fc$ is supported by a fiber of $\pi$ which is a complex line. A positive $\ddc$-closed current on a Riemann surface is defined by a positive harmonic function. But on a complex line, the only positive harmonic functions are the constants. We conclude that extremal elements in $\Fc$ are multiples of currents of integration on fibers of $\pi$. In order to get the lemma, we only need to show that $S$ is an average of those extremal currents.

Observe that it is enough to consider the local setting. We will consider the following slightly more general situation where $\pi$ is the projection 
$(z_1,z_2)\mapsto z_1$ from $\C\times W$ to $\C$ and $S$ is a vertical positive $\ddc$-closed current with support in $\overline\D\times W$. Here $W$ is any connected open set in $\C$. For our lemma, we only need the case $W=\C$. Fix also a closed disc $\overline D$ in $W$.

The mass of $S$ on a Borel set $K$ is the mass of the measure $S\wedge (idz_1\wedge d\overline z_1+idz_2\wedge d\overline z_2)$ on $K$. The last measure is equal to $S\wedge idz_2\wedge d\overline z_2$ because $S$ is vertical.
On the other hand, if $\pi'(z_1,z_2):=z_2$ is the projection from $\C\times W$ to $W$, $\pi'_*(S)$ is a positive $\ddc$-closed $(0,0)$-current on $W$. So it is defined by a positive harmonic function $h$ on $W$. We see that the mass of $S$ on $\C\times U$ is equal to the integral of $h (i dz_2\wedge d\overline z_2)$ on $U$ for every relatively compact open subset $U$ of $W$. Thus, by Harnack's inequality, applied for $h$, this mass is bounded by the mass of $S$ on $\C\times D$ times a constant which depends only on $U$. 

Finally, consider the convex cone of positive $\ddc$-closed vertical currents $S$ as above. 
The last discussion on the mass of $S$ implies that the set of currents with mass 1 on $\C\times D$ is compact and is a basis of the considered cone. 
Therefore, Choquet's representation theorem implies that any current in the cone is an average on the extremal elements. The lemma follows.
\endproof

\medskip
\noindent
{\bf Proof of Proposition \ref{prop_tangent_local}.} Consider a smooth test 2-form $\Omega$ with compact support in $\D\times\C$. Denote by $f(z)$ the coefficient of $dz_2\wedge d\overline z_2$ in $\Omega$.  
By definition of $S$ and the above discussion on negligible families of forms, we have 
$$\langle S,\Omega\rangle = \lim_{n\to\infty} \big\langle T,\tau^*A_{\lambda_n}^*(f(z)dz_2\wedge d\overline z_2)\big\rangle = \lim_{n\to\infty}\big\langle T, |\lambda_n|^2f(\tau_1,\lambda_n\tau_2) d\tau_2\wedge d\overline \tau_2 \big\rangle$$
and
$$\lim_{n\to\infty} \big\langle T,a_{\lambda_n}^*\Omega\big\rangle = \lim_{n\to\infty} \big\langle T,a_{\lambda_n}^*(f(z)dz_2\wedge d\overline z_2)\big\rangle = \lim_{n\to\infty}\big\langle T, |\lambda_n|^2f(z_1,\lambda_nz_2) dz_2\wedge d\overline z_2 \big\rangle.$$
So it is enough to check that the family
$$|\lambda|^2f(\tau_1,\lambda\tau_2) d\tau_2\wedge d\overline \tau_2 - |\lambda|^2f(z_1,\lambda z_2) dz_2\wedge d\overline z_2$$
is negligible. But this can be easily deduced from the local description of $\tau$ and $d\tau$ and the fact that $|z_2|\lesssim |\lambda|^{-1}$ in the support of the above form.
\hfill $\square$

\medskip

In the rest of this section, we complete the proof of Theorem \ref{th_tangent}. It remains to compute the mass of the measure $\nu$. 
As mentioned at the end of the last section, the above tangent current $S$ does not depend on the choice of $\tau$. From now on, the map $\tau$ is the one constructed in Example \ref{ex_admissible}. The following result shows that it is closer to holomorphic maps than general admissible maps. 

\begin{proposition} \label{prop_ex_admissible}
In local coordinates $z=(z_1,z_2)$ as above, we have for $z_2\to 0$
$$\tau(z_1,z_2)= \Big(z_1+ O(|z_2|^2), z_2+ O(1)z_2^2+O(|z_2|^3)\Big)$$
and
 $$d\tau(z_1,z_2)= \Big(dz_1+O^*(|z_2|^2), dz_2+O^*(|z_2|^2)\Big).$$
 \end{proposition}
\proof
The second assertion is a consequence of the first one. Observe also that the first identity is equivalent to the similar identity for $\widetilde\tau:=\tau^{-1}$. We will prove the last one.
Since $d\widetilde\tau(z_1,z_2)$ is the identity when $z_2=0$, we have
$\widetilde\tau=\id + O(|z_2|^2)$. So if we write $\widetilde\tau=(\widetilde\tau_1,\widetilde\tau_2)$ in coordinates $(z_1,z_2)$, we only have to check that $\widetilde\tau_2=z_2+ O(1)z_2^2+O(|z_2|^3)$. This property means there are no terms with $z_2\overline z_2$ or $\overline z_2^2$ in the Taylor expansion of $\widetilde\tau_2$ in $z_2,\overline z_2$ with functions in $z_1$ as coefficients. So it is enough to check it on each complex line $\{z_1\}\times\C$. 
 Recall that in the local coordinates $(z_1,z_2)$ as above, we identify this complex line with the fiber of $N_{V|X}$ over $(z_1,0)$.
We will need to make some changes of coordinates. So we first check that the property does not depend on our choice of coordinates. 

Now consider another system of local holomorphic coordinates $(z_1',z_2')$ such that $z_2'=0$ on $V$. We can write $z_1'=H(z_1,z_2)$ and $z_2'=\beta(z_1) z_2+h(z_1,z_2)z_2^2$, where $H,\beta$ and $h$ are holomorphic functions. 
For $a'=H(a,0)$, the two complex lines $\{a\}\times\C$ for the coordinates $(z_1,z_2)$ and $\{a'\}\times \C$ for the coordinates $(z_1',z_2')$ are both  identified with the same fiber of $N_{V|X}$. The linear map relying them is  $(a,z_2)\mapsto (a',\beta(a) z_2)$. We will keep the notation $\widetilde\tau=(\widetilde\tau_1,\widetilde\tau_2)$ for the map $\widetilde\tau$ in coordinates $(z_1,z_2)$ and use $\widetilde\tau'=(\widetilde\tau_1',\widetilde\tau_2')$ for the same map in coordinates $(z_1',z_2')$. With these notations, the point
$\widetilde\tau(a,\beta(a)^{-1}b')$ in  coordinates $(z_1,z_2)$ and the point $\widetilde\tau'(a',b')$ in  coordinates $(z_1',z_2')$ represent the same point of $X$. It follows that 
$$\widetilde\tau'_2(a',b')=\beta(a)\widetilde\tau_2(a,\beta(a)^{-1}b') +h\big(\widetilde\tau(a,\beta(a)^{-1}b')\big)\widetilde\tau_2(a,\beta(a)^{-1}b')^2.$$
We see that if $\widetilde\tau_2(a,b)=b+O(1)b^2+O(|b|^3)$ then $\widetilde\tau_2'(a',b')$ satisfies a similar property.

In the rest of the proof, we show that
 $\widetilde\tau_2=z_2+ O(1)z_2^2+O(|z_2|^3)$.
Without loss of generality, we will only check the property for $z_1=0$ and $z_2=tw$ with $t\in\R_+$ and $|w|=1$. 
We need that $\tau$ is associated with a K\"ahler form $\omega$ as explained in Example \ref{ex_admissible}.
 
In a neighbourhood of 0, we can write $\omega=i\ddbar u$ with $u$ a smooth strictly psh function. Subtracting from $u$ a pluriharmonic function, we can assume the existence of a positive definite $2\times 2$-matrix  $(c_{ij})$ such that 
$$u(z)=c_{11}|z_1|^2+c_{12}z_1\overline z_2+c_{21} \overline z_1 z_2 +c_{22} |z_2|^2+O(\|z\|^3).$$ 
We will make changes of coordinates keeping the property $V=\{z_2=0\}$.  With a linear change of coordinates 
$(z_1,z_2)\mapsto (\gamma_1z_1+\gamma_2z_2,\beta z_2)$, we can assume that 
$$u(z)=|z_1|^2+ |z_2|^2+O(\|z\|^3).$$ 
Then, using a change of coordinates of type $z_2\mapsto z_2+\gamma z_1z_2+\beta z_2^2$, we can assume that the coefficients of $z_1z_2\overline z_2$,  $z_2^2\overline z_2$ and their conjugates in the last $O(\|z\|^3)$ vanish. Note that since $u$ is real, when we  eliminate the coefficient of a monomial, the coefficient of its complex conjugate is also eliminated. 
Next, using a change of coordinates of type $z_1\mapsto z_1+\gamma z_1^2+\beta z_1z_2+\theta z_2^2$, we can assume that the coefficients of $z_1^2\overline z_1$, $z_1z_2\overline z_1$, $z_2^2\overline z_1$ and their conjugates vanish.
It follows that 
$$\omega=idz_1\wedge d\overline z_1 + idz_2\wedge d\overline z_2 + O(|z_1|) dz_1\wedge d\overline z_2 + O(|z_1|) dz_2\wedge d\overline z_1 +O(\|z\|^2).$$

For the rest of the proof, we use real coordinates $x=(x^1,x^2,x^3,x^4)$ such that $z_1=x^1+ix^2$ and $z_2=x^3+ix^4$. Denote by $v=(v^1,v^2,v^3,v^4)$ the unit tangent vector to $X$ at 0 corresponding to $(0,w)\in\C^2$, i.e. $v^1=v^2=0$ and $w=v^3+iv^4$. 
So $\widetilde\tau(0,tw)$ is equal to $\exp(tw)$, where $\exp$ denotes the exponential map from the tangent space to $X$ at 0 to $X$. 
If we write $\widetilde \tau(0,tw)=(x^1(t),x^2(t),x^3(t),x^4(t))$, then $x^j(t)$  satisfy the geodesic equations
$$\ddot x^j= -\Gamma^j_{kl} \dot x^k \dot x^l \qquad \text{and} \qquad \dot x^j(0)=v^j,$$
where $\Gamma^j_{kl}$ are the Christoffel symbols associated with the considered K\"ahler metric.

We will show in the present setting that $\widetilde\tau_2(0,tw)=tw+ O(t^3)$ and we already know that $\widetilde\tau_2(0,tw)=tw+ O(t^2)$. This is equivalent to checking that $\ddot x^j(0)=0$ for $j=3,4$. Note that the property implies that there is no term of order 2 in the Taylor expansion of $\widetilde\tau_2(0,z_2)$ in the latest system of coordinates. According to the discussion at the beginning of the proof, a term with $z_2^2$ may appear when we come back to the original coordinates. 
Since $v^j=0$ for $j=1,2$, we only need to show that $\Gamma^j_{kl}(0)=0$ for $j,k,l\in\{3,4\}$. 
Let  $g=(g_{jk})$ be
the Riemannian metric associated with $\omega$. The above description of $\omega$ implies that 
$g_{jk}=\delta_{jk}+O(|x^1|+|x^2|+\|x\|^2)$ for all $j,k$, where $\delta_{jk}=1$ if $j=k$ and 0 otherwise. The coefficients of the inverse $(g^{jk})$ of the matrix $(g_{jk})$ satisfy a similar property. 
Recall that the Christoffel symbols are given by 
$$\Gamma^j_{kl}={1\over 2} g^{jm} \Big({\partial g_{mk}\over \partial x^l} + {\partial g_{ml}\over \partial x^k} -{\partial g_{kl}\over \partial x^m}\Big).$$
It is now easy to check that $\Gamma^j_{kl}(0)=0$ for $j,k,l\in\{3,4\}$. The proposition follows. 
\endproof

\noindent
{\bf End of the proof of Theorem \ref{th_tangent}.} It remains to compute the mass of $\nu$. 
Let $\Phi$ be any smooth closed $(1,1)$-form on a small neighbourhood of $V$ in $N_{V|X}$ which is cohomologous to $\{V\}$. It is enough to check that $\langle T_\lambda,\Phi\rangle$ converges to $\{T\}\smallsmile \{V\}$ because with the above notation $\langle T_{\lambda_n},\Phi\rangle$ tends to the mass of $\nu$. We first construct a suitable form $\Phi$.

Let $\phi$ be a function with support in a small neighbourhood of $V$ in $N_{V|X}$ such that in local coordinates as above $\phi-\log |z_2|$ is smooth. This property does not depend on the choice of local coordinates $(z_1,z_2)$ with $z_2=0$ on $V$. So it is easy to construct such a function $\phi$ using a partition of unity on a neighbourhood of $V$. 
Define $\Phi:=-\ddc \phi+[V]$. This is a smooth $(1,1)$-form in a neighbourhood of $V$ which is cohomologous to $\{V\}$. 

Define $\phi_\lambda:=\phi\circ A_\lambda$ and $\psi_\lambda:=\phi_\lambda\circ \tau$. Observe that in local coordinates, we also have 
that $\phi_\lambda-\log|z_2|$ is smooth and hence $\psi_\lambda-\log|\tau_2|$ is also smooth. Define
$\Phi_\lambda:=(A_\lambda)^*(\Phi)=-\ddc \phi_\lambda+[V]$, 
  $\Psi_\lambda:=\tau^*(\Phi_\lambda)$ and $\Psi'_\lambda:=-\ddc \psi_\lambda+[V]$. We have to show that $\langle T,\Psi_\lambda\rangle\to \{T\}\smallsmile \{V\}$. Proposition \ref{prop_ex_admissible} implies that 
 $$\log|\tau_2|-\log|z_2|=\log \Big|1+O(1)z_2 +{O(|z_2|^3)\over z_2}\Big|.$$
 Since the function $\xi\mapsto\log |\xi|$ is smooth near the point $1\in\C$ and the function $O(|z_2|^3)/ z_2$ is of class $\Cc^{1,1}$, the function
   $\log|\tau_2|-\log|z_2|$ is also $\Cc^{1,1}$. Therefore,  $\Psi_\lambda'$ is a closed $(1,1)$-current with $\Cc^{1,1}$ potentials. By Lemma \ref{lemma_harm_c11}, we have 
 $$\langle T, \Psi'_\lambda\rangle_{X\setminus V}=\{T\}\smallsmile \{\Psi_\lambda'\}=\{T\}\smallsmile \{V\}.$$ 
Recall that $T$ has no mass on $V$ and $\Psi_\lambda,\Psi_\lambda'$ are $L^\infty$ forms. So it is enough to check that $\langle T,\Psi_\lambda-\Psi'_\lambda\rangle_{X\setminus V}\to 0$. 

For that purpose, we will use local coordinates and apply Lemma \ref{lemma_negligible}. 
We want to show that $\Psi_\lambda-\Psi'_\lambda$ is negligible in the sense of Definition \ref{def_negligible}. The condition on the support of this form is clearly satisfied. 
We will only consider the coefficients of $idz_2\wedge d\overline z_2$ because the other coefficients can be treated in the same way.
Define $\widetilde\phi:=\phi-\log|z_2|$. This is a smooth function. We have outside $V$
$$\Psi_\lambda-\Psi'_\lambda = -\tau^*\ddc (\widetilde\phi\circ A_\lambda) + \ddc (\widetilde \phi\circ A_\lambda\circ \tau)+\ddc \log|\tau_2|.$$
Recall that $\log |z_2|-\log|\tau_2|$ is $\Cc^{1,1}$. Therefore, the coefficient of $idz_2\wedge d\overline z_2$ in $\ddc \log|\tau_2|$ is bounded and hence negligible.  

The map $\tau$ is described in
Proposition \ref{prop_ex_admissible}. A direct computation shows that, modulo a negligible family, the coefficient of $idz_2\wedge d\overline z_2$  in $\Psi_\lambda-\Psi'_\lambda$
is equal to $1/\pi$ times 
$$|\lambda|^2{\partial \widetilde\phi \over \partial z_2\partial \overline z_2}(\tau_1,\lambda \tau_2)-{\partial^2\over \partial z_2\partial \overline z_2} \big[\widetilde\phi(\tau_1,\lambda \tau_2)\big].$$
Using again the local description of $\tau$, we see that the last expression is negligible.
\hfill $\square$

\medskip

Note that we need Proposition \ref{prop_ex_admissible} in order to get the $\Cc^{1,1}$ property of $\log|\tau_2|-\log|z_2|$. 
If there exists a {\bf holomorphic} admissible map, we can use it instead of the above map $\tau$. This is the case when $X=\P^2$ and $V$ is a projective line in $X$.

\section{Directed currents near a singularity of a foliation} \label{section_sing}

Let $\Fc$ be a holomorphic foliation on a compact K\"ahler surface. Let $V$ be a smooth compact complex curve in $X$ which is invariant by $\Fc$. 
Assume that all singular points of $\Fc$ in $V$, when they exist, are hyperbolic. Consider also a positive $\ddc$-closed $(1,1)$-current $T$ directed by $\Fc$. 
Assume that $T$ has no mass on $V$ and on the separatrices of $\Fc$ at every singular point $p\in V$, see also Propositions \ref{prop_T_V} and \ref{prop_current_leaf}.
In the next section, we will show that the tangent currents to $T$ along $V$ vanish. This result requires properties of $T$ near these singular points that we will discuss below.

Let $p\in V$ be a singular point of the foliation. Since $p$ is hyperbolic, there is a local coordinate system $z=(z_1,z_2)$ centered at $p$, with $|z_1|<3$, $|z_2|<3$,  such that in the bidisc $\D(0,3)^2=\{|z_1|<3,|z_2|<3\}$,  the foliation $\Fc$ is defined by the vector field
$$F:=z_2{\partial \over \partial z_2}-\eta z_1{\partial\over \partial z_1} \ \raisepunct{,}$$
where $\eta=a+ib$, $a,b\in\R$, and $b\not=0$.
Observe that the two axes of the bidisc $\D(0,3)^2$ are invariant and are the separatrices of the foliation in the bidisc $\D(0,3)^2$. 
So $V$ is one of these separatrices. We can assume that $V$ is given in $\D(0,3)^2$ by the equation $z_2=0$.

Consider the ring $\A$ defined  by
$$\A:=\big\{\alpha\in\C, \ e^{-2\pi|b|}<|\alpha|\leq 1\big\}.$$
Define also the sectors $\S$ and $\S'$  by
$$\S:=\big\{\zeta=u+iv\in\C,\quad  v>0 \quad\text{and}\quad bu+av >0 \big\}$$
and
 $$\S':=\big\{\zeta=u+iv\in\C,\quad   v>-\log 3 \quad\text{and}\quad bu+av >-\log3 \big\}\supset \S.$$
Note that the sector $\S$ is contained in the upper half-plane $\H:=\{u+iv,\ v>0\}$. 
For $\alpha\in \C^*$, consider the following manifold $\L_\alpha$ immersed in $\C^2$ and defined by 
$$z_1=\alpha e^{i\eta(\zeta+\log|\alpha|/b)} \quad \text{and} \quad z_2=e^{i(\zeta+\log|\alpha|/b)} \quad \text{with} \quad \zeta=u+iv\in\C.$$
Note that the map $\zeta\mapsto (z_1,z_2)$ is injective because $\eta\not\in\R$.
It is easy to check the following properties  
\begin{enumerate}
\item[(1)] $\L_\alpha$ is tangent to the vector field $F$ and is a submanifold of $\C^{*2}$.
\item[(2)] $\L_{\alpha_1}$ is equal to $\L_{\alpha_2}$ if $\alpha_1/\alpha_2=e^{2ki\eta\pi}$ for some $k\in\Z$ and they are disjoint otherwise. In particular, $\L_{\alpha_1}$ and $\L_{\alpha_2}$ are disjoint if $\alpha_1,\alpha_2\in\A$ and $\alpha_1\not =\alpha_2$.
\item[(3)] The union of $\L_\alpha$ is equal to $\C^{*2}$ for $\alpha\in \C^*$, and then also for $\alpha\in \A$.
\item[(4)] The intersection $L_\alpha:=\L_\alpha\cap \D^2$ of $\L_\alpha$ with the unit bidisc $\D^2$ is given by the same equations as in the definition of $\L_\alpha$ but with $\zeta\in \S$. Moreover, $L_\alpha$ is a connected submanifold of  $\D^{*2}$. In particular,  it is a leaf of $\Fc\cap\D^2$.
\item[(5)] Similarly, the intersection $L_\alpha':=\L_\alpha\cap \D(0,3)^2$  is given by the same equations with $\zeta\in \S'$. Moreover, $L_\alpha'$ is a connected submanifold of  $\D^*(0,3)^{2}$ and is the leaf of $\Fc\cap\D(0,3)^2$ which contains $L_\alpha$.
 \end{enumerate}
 
From now on, assume for simplicity that $b>0$ and hence $\S$ is a sector of angle  $\arctan (-a/b)\in (0,\pi)$ having $\R_+$ in its boundary (otherwise, it has $\R_-$ in the boundary and we need to change slightly the definition of the map $\Phi$ below). 
Consider the map $\Phi:\S\to\H$ with 
$$\Phi(\zeta):=\zeta^\gamma \quad \text{with} \quad \gamma:={\pi\over \arctan(-b/a)}>1,$$
sending $\S$ bi-holomorphically to the upper half-plane $\H\subset\C$.

\begin{lemma} \label{lemma_decom_D2}
There is a positive measure $\mu$ of finite mass on $\A$ and positive harmonic functions $h_\alpha$ on $L_\alpha'$ for $\mu$-almost every $\alpha\in\A$ such that in $\D(0,3)^2$
$$T=\int_\A T_\alpha d\mu(\alpha), \quad \text{where} \quad T_\alpha:=h_\alpha[L_\alpha'].$$
Moreover, the mass of $T_\alpha$ in $\D(0,2)^2$ is $1$ for $\mu$-almost every $\alpha\in\A$.
\end{lemma}
\proof
Recall that $T$ has no mass on the separatrices. 
Let $\Sigma$ denote the quotient of $\C^*$ by the equivalence relation $\alpha_1\sim\alpha_2 \Longleftrightarrow \alpha_1/\alpha_2=e^{ki\eta\pi}$ for some $k\in\Z$. This is a complex torus of dimension 1 and the natural projection from $\A$ to $\Sigma$ is bijective.
 
According to the above description of $\Fc$ in $\D(0,3)^2$ the map $\pi$ which associates to a point in $L_\alpha'$  the image of $\alpha$ in $\Sigma$ is a holomorphic map and $T$ is directed by the fibers of $\pi$.  It follows that if $\Omega$ is any smooth form of degree 1 or 2 on $\Sigma$, then $T\wedge \pi^*(\Omega)=0$. Arguing as in  Lemma \ref{lemma_decomp},
we obtain the decomposition of $T$ as stated in the lemma. Finally, we multiply $\mu$ by a suitable 
positive function $\theta(\alpha)$ and divide $h_\alpha$ by $\theta(\alpha)$ in order to assume that the mass of $T_\alpha$ 
in $\D(0,2)^2$ is 1 for $\mu$-almost every $\alpha\in\A$. 
\endproof

Define 
$$H_\alpha(\zeta):=h_\alpha\Big(\alpha e^{i\eta(\zeta+\log|\alpha|/b)}, e^{i(\zeta+\log|\alpha|/b)}\Big).$$
This is a positive harmonic function in $\S'$.
Define also $\widetilde H_\alpha:=H_\alpha\circ \Phi^{-1}$. This is a positive harmonic function in $\H$ which is continuous up to the boundary because $\S'\supset \overline\S=\Phi^{-1}(\overline\H)$.
The following lemma was obtained in \cite{FS2}. We will give a slightly different proof, see also \cite{Nguyen1}. 

\begin{lemma} \label{lemma_poisson}
For $\mu$-almost every $\alpha\in\A$, the function $\widetilde H_\alpha$ is the Poisson integral of its boundary values, that is, for $U,V\in\R$ with $V>0$, we have
$$\widetilde H_\alpha(U+iV)={1\over\pi}\int_{x\in\R} \widetilde H_\alpha(x){V\over V^2+(x-U)^2} dx.$$
\end{lemma}
\proof
We only consider generic values of $\alpha$ with respect to $\mu$. In particular, the mass of $T_\alpha$ on $\D(0,2)^2$ is 1. 
By definition, the mass of $T_\alpha$ in $\D^2$ is the mass of the following positive measure in $\D^2$:
$$T_\alpha\wedge (idz_1\wedge d\overline z_1+idz_2\wedge d\overline z_2) = h_\alpha(z) (idz_1\wedge d\overline z_1+idz_2\wedge d\overline z_2)\wedge [L_\alpha].$$
Using the above parametrization of $L_\alpha$ by $\S$, we see that the mass of this measure is equal to the one of it pull-back to $\S$ which is equal on $\S$ to the measure
$$H_\alpha(\zeta) \big(|\eta|^2| e^{-2(bu+av)} +e^{-2v}\big)id\zeta\wedge d\overline\zeta, \quad \text{where we write}Ê\quad \zeta=u+iv.$$

So the considered mass of $T_\alpha$ is the integral of the last expression on $\S$ and it is larger than or equal to the integral on the sub-domain
$$\S_{1}:=\{1/2<v<1, u>M\}, \quad \text{where }  M>0 \text{ is a fixed large constant}.$$
Hence, it is bounded below by
$$e^{-2}\int_{\S_1} H_\alpha(\zeta)id\zeta\wedge d\overline\zeta$$
since $|\eta|^2 e^{-2(bu+av)} +e^{-2v}\geq e^{-2}$ on $\S_1$. So the last integral is bounded by a constant.

Now, using the change of variables $Z=U+iV:=\Phi(\zeta)=\zeta^\gamma$ and $\zeta=\Phi^{-1}(Z)=Z^{1/\gamma}$, 
we have $\widetilde H_\alpha(U+iV)=H_\alpha(u+iv)$. Since $\widetilde H_\alpha$ is a positive harmonic function on $\H$ and continuous up to the boundary of $\H$, it has the following form 
$$\widetilde H_\alpha(U+iV)=c_\alpha V + {1\over\pi}\int_{x\in\R} \widetilde H_\alpha(x){V\over V^2+(x-U)^2} dx$$
for some constant $c_\alpha\geq 0$, see \cite[Th. 16.1.7]{Hormander}. It remains to check that $c_\alpha=0$.
It is not difficult to see that for $\zeta\in \S_1$, we have $V\gtrsim  u^{\gamma-1}$, see Lemma \ref{lemma_Phi_Lambda_infty} below. Hence, 
$$\int_{\S_1} H_\alpha(\zeta)id\zeta\wedge d\overline\zeta\gtrsim \int_{\S_1} c_\alpha u^{\gamma-1} du dv.$$
Since the first integral is finite, we deduce that  $c_\alpha=0$.
\endproof

\begin{lemma} \label{lemma_poisson_est}
We have for $\mu$-almost every $\alpha\in\A$
$$\int_{x\in\R} \widetilde H_\alpha(x) |x|^{-1+1/\gamma} dx  \leq C,$$
where $C>0$ is a constant independent of  $\alpha$.  
\end{lemma}
\proof
We use the notation introduced in the previous lemma.
Since $H_\alpha$ is harmonic and positive on $\S'$, by Harnack's inequality, if $\zeta=u+iv\in \S_1$, the ratio $H_\alpha(\zeta)/H_\alpha(u-M)$ is bounded below by a strictly positive constant independent of $\alpha$. Recall that the integral of $H_\alpha$ on $ \S_1$ is bounded by a constant. We then deduce that its integral on $\R_+$ is also bounded by a constant. Using the relation $\widetilde H_\alpha(x)=H_\alpha(x^{1/\gamma})$, we obtain that the integral
$$\int_0^\infty \widetilde H_\alpha(x) |x|^{-1+1/\gamma} dx$$
is bounded by a constant.

Similarly, we obtain the same property for the integral on $(-\infty,0]$ by using the half-line $\{au+bv=0, v\geq 0\}$ on the boundary of $\S$ instead of the half-line $\R_+$ used above. It is sent by $\Phi$ to $\R_-$. The lemma follows.
\endproof

\begin{lemma} \label{lemma_H_bounded}
For $\mu$-almost every $\alpha\in\A$, the function $\widetilde H_\alpha$ is bounded on $\R$ by a constant $C>0$ which is independent of $\alpha$.
\end{lemma}
\proof
It is enough to show that $H_\alpha$ is bounded by a constant (independent of $\alpha$) on the boundary of $\S$. For $u\in\R_+$, by Harnack's inequality, $H_\alpha(u+t)$ is bounded by a constant times $H_\alpha(u)$ when $0\leq t\leq 1$. We have seen that the integral of $H_\alpha$  on $\R_+$ is bounded by a constant. It is then easy to deduce that $H_\alpha(u)$ is bounded by a constant. Similarly, $H_\alpha$ is also bounded on the other edge of the cone $\S$. The lemma follows.
\endproof

\begin{lemma} \label{lemma_s_n}
For each integer $n\geq 1$, there is a  real number $s_n\geq n$ such that  $J_n(s_n)\leq n^{-1}s_n^{\gamma-1}/\log s_n$, where
$$J_n(s):=\int_{\alpha\in \A} \int_{|x|\leq ns^\gamma}\widetilde H_\alpha(x)dx d\mu(\alpha).$$ 
\end{lemma}
\proof
Assume that the estimate does not hold. Then there is $n\geq 1$ such that $J_n(s)\geq n^{-1} s^{\gamma-1}/\log s$ for $s\geq n$. 
A direct computation and Fubini's theorem imply that
\begin{eqnarray*}
\int_{s\geq n} J_n(s) s^{-\gamma} ds  & \leq &  \int_{s\geq n}\int_{\alpha\in \A} \int_{|x|\leq ns^\gamma} \widetilde H_\alpha(x)s^{-\gamma}dx d\mu(\alpha)ds \\
 & \leq &  \int_{\alpha\in \A} \int_{x\in \R} \int_{s\geq |x/n|^{1/\gamma}} \widetilde H_\alpha(x)s^{-\gamma}dx d\mu(\alpha)ds \\
& \lesssim  &  \int_{\alpha\in \A} \int_{x\in\R}\widetilde H_\alpha(x)|x|^{-1+1/\gamma}dx d\mu(\alpha).
\end{eqnarray*}
By Lemma \ref{lemma_poisson_est}, the last double integral is finite. 
On the other hand, since $J_n(s) s^{-\gamma}\geq n^{-1}/(s\log s)$, the first integral is bounded below by 
$$n^{-1}\int_{s\geq n} {ds\over s\log s} =n^{-1}\log\log s\Big|_{s=n}^\infty =\infty.$$
This is a contradiction. The lemma follows.

Note that the same arguments show that for every $\epsilon>0$ the set 
$$E:=\{s\geq 2, \quad J_n(s)\geq \epsilon s^{\gamma-1}/\log s\}$$
has finite length with respect to the measure $ds/(s\log s)$. An analogous property is well-known in the Nevanlinna theory. 
\endproof

\begin{remark} \rm
When $V$ contains more than 1 singular point of the foliation, denote by $p_1,\ldots,p_N$ these singular points. 
Denote also by $\gamma_1,\ldots,\gamma_N$ the corresponding angles $\gamma$ associated with $p_1,\ldots,p_N$ respectively. 
Define for each singularity $p_i$ an integral $J_n^{(i)}(s)$ as in the last lemma. It follows from the proof of Lemma \ref{lemma_s_n} that there exists $s_n\geq n$ such that $\sum_i J_n^{(i)}(s_n)s_n^{-\gamma_i}\leq n^{-1}/(s_n\log s_n)$. So we can choose the same $s_n$ satisfying the above lemma for each $p_i$. 
\end{remark}

Let $\rho>0$ and $N\geq 2\rho$ be constants. We will use the following estimate in the next section for a particular value of $\rho$ depending only on $\gamma$ and for  a constant $N$ large enough.

\begin{lemma} \label{lemma_s_n_bis}
Let $s_n$ be as in Lemma \ref{lemma_s_n}. Then $K(s_n)=o(s_n^{\gamma-1})$ as $n\to\infty$, where
$$K(s):=\int_{\alpha\in \A} \int_{|x|\leq Ns^\gamma} \widetilde H_\alpha(x)\log{6Ns^\gamma\over |x+\rho s^\gamma |}dx d\mu(\alpha).$$
\end{lemma}
\proof
We only consider $n$ large enough such that $n\geq N$. So $s_n$ is also a large number.
Observe that the function under the considered integral is positive when $s>0$. Define
$$\widetilde H(x):=\int_{\alpha\in \A} \widetilde H_\alpha(x)d\mu(\alpha).$$
By Lemma \ref{lemma_H_bounded}, this function is bounded by a constant $C$. We also have
$$J_n'(s_n):= \int_{|x|\leq Ns_n^\gamma}\widetilde H(x)dx\leq  \int_{|x|\leq ns_n^\gamma}\widetilde H(x)dx=J_n(s_n)$$
and 
$$K(s_n) = \int_{|x|\leq N s_n^{\gamma}}\widetilde H(x)\log{6Ns_n^\gamma\over |x+\rho s_n^\gamma|}dx.$$

Denote by $\widetilde H'(x)$ the function which is equal to $C$ times the characteristic function on the interval centered at $-\rho s_n^\gamma$ and of length $2l_n:=J'_n(s_n)/C$. Note that  by Lemma \ref{lemma_s_n}, 
this interval is contained in the interval $\{|x|\leq Ns_n^\gamma\}$ where we consider the integrals. 
Let $K'(s_n)$ be the integral similar to the last expression of $K(s_n)$ but we replace $\widetilde H$ with $\widetilde H'$. Observe that 
$$\widetilde H'(x)- \widetilde H(x)\geq 0 \quad \text{and} \quad  \log{6Ns_n^\gamma\over |x+\rho s_n^\gamma|} \geq \log{6Ns_n^\gamma\over l_n} \quad \text{in} \quad  [-\rho s_n^\gamma-l_n,-\rho s_n^\gamma+l_n]$$ 
and  
$$\widetilde H'(x)- \widetilde H(x)\leq 0 \quad \text{and} \quad  \log{6Ns_n^\gamma\over |x+\rho s_n^\gamma|} \leq \log{6Ns_n^\gamma\over l_n} \quad \text{outside} \quad  [-\rho s_n^\gamma-l_n,-\rho s_n^\gamma+l_n].$$ 
In both cases, we get
$$(\widetilde H'(x)- \widetilde H(x))\log{6Ns_n^\gamma\over |x+\rho s_n^\gamma|}\geq (\widetilde H'(x)- \widetilde H(x))\log{6Ns_n^\gamma\over l_n}\cdot$$
It follows that
$$K'(s_n)-K(s_n) \geq  \int_{|x|\leq N s_n^{\gamma}}(\widetilde H'(x)-\widetilde H(x))\log{6Ns_n^\gamma\over l_n}dx.$$
Using the definition of $l_n$, we obtain that the last integral is equal to  0. Therefore, we have $K'(s_n)\geq K(s_n)$.

By Lemma \ref{lemma_s_n}, we have $J'_n(s_n)=o(s_n^{\gamma-1}/\log s_n)$. Define  
$$\delta_n:=J'_n(s_n)/(12CNs_n^{\gamma})=o(1/(s_n\log s_n)).$$
Hence
$$K(s_n)\leq K'(s_n) = C\int_{|x+\rho s_n^\gamma|\leq \delta_n 6Ns_n^\gamma}
\log{6Ns_n^\gamma\over |x+\rho s_n^\gamma|}dx.$$
Using the variable $t:=(x+\rho s_n^\gamma)/(6Ns_n^\gamma)$, we obtain that the last integral is equal to a constant times
$$s_n^\gamma\int_{|t|\leq \delta_n} -\log |t| dt
\lesssim -s_n^\gamma\delta_n\log\delta_n =o(s_n^{\gamma-1}).$$
The lemma follows.
\endproof

\section{Unique ergodicity for foliations with invariant curves}

In this section, we will prove the following theorem which implies the results stated in the introduction.

\begin{theorem} \label{th_main_bis}
Let $\Fc$ be a holomorphic foliation by Riemann surfaces in a compact K\"ahler surface $X$.  Let $V$ be a compact complex curve in $X$ invariant by $\Fc$ and $T$ a positive $\ddc$-closed $(1,1)$-current directed by $\Fc$ and having no mass on $V$. If $\Fc$ has singularities on $V$, we assume they  are all hyperbolic and that $T$ has no mass on their separatrices. Then we have 
$\{T\}\smallsmile \{V\}=0$.
\end{theorem}

Recall that a theorem by Jouanolou implies that the number of invariant compact curves is finite provided that the foliation admits no meromorphic first integral \cite{Jouanolou2}. This is the case when the foliation contains a hyperbolic singularity. 

\begin{corollary} \label{cor_main_1}
Let $X$ and $\Fc$ be as in Theorem \ref{th_main_bis} with only hyperbolic singularities. Let $E$ denote the vector subspace of $H^{1,1}(X,\R)$ spanned by the classes of compact complex curves invariant by $\Fc$. Assume that $E$ contains the class of a continuous real closed $(1,1)$-form $\alpha$ on $X$ which is strictly positive along the leaves of $\Fc$, e.g. $\alpha$ is a K\"ahler form. Then every positive $\ddc$-closed $(1,1)$-current directed by $\Fc$ is a linear combination of currents of integration on irreducible invariant curves (there are finitely many such curves).  
\end{corollary}
\proof
Let $V_1,\ldots, V_n$ be the invariant  irreducible compact curves. So $\{\alpha\}$ is a linear combination of $\{V_1\},\ldots,\{V_n\}$. Let $T$ be a positive $\ddc$-closed $(1,1)$-current directed by $\Fc$. 
We  show that $T$ is a linear combination of currents of integration on invariant compact curves.
We can assume that $T$ has no mass on $V_i$ and on the separatrices, see Proposition \ref{prop_current_leaf}. By Theorem \ref{th_main_bis}, we have $\{T\}\smallsmile \{V_i\}=0$. By Theorem \ref{th_tangent}, we get $\{T\}\smallsmile \{\alpha\}=0$. Since $\alpha$ is continuous and closed, we deduce that $\{T\wedge\alpha\}=0$. But $\alpha$ is strictly positive along the leaves of $\Fc$ and $T$ is positive directed by $\Fc$, hence $T=0$.
\endproof

 \noindent
{\bf Proof of Theorem \ref{th_main_P2_bis}.}
Observe that $\dim_\R H^{1,1}(\P^2,\R)=1$ and hence the class of a complex compact curve in $\P^2$ is always K\"ahler.  
If $\Fc$ admits an invariant algebraic curve, we obtain (2) as a consequence of Corollary \ref{cor_main_1}. Otherwise, $\Fc$ has no invariant algebraic curve and property (1) is exactly the main theorem in \cite{FS2}. 
\hfill $\square$ 
 
 \medskip

 \noindent
{\bf Proof of Theorem \ref{th_main_P2}.} The first assertion is a consequence of Theorem  \ref{th_main_P2_bis}. The second assertion follows from the fact that any limit of $\tau_r^w$ is a positive current of mass 1 directed by $\Fc$, which is $\ddc$-closed thanks to Proposition \ref{prop_current_Nev}. 
\hfill $\square$ 
 
 \medskip
 
In the rest of this section, we give the proof of Theorem \ref{th_main_bis}. 
We first consider the case where $V$ is smooth. By Theorem \ref{th_tangent}, we only need to show that the tangent current to $T$ along $V$ is zero.
Consider the sequence $\lambda_n:=e^{s_n}$ with  $s_n$ defined  in Section \ref{section_sing}. By extracting a subsequence, we can assume that there is  a tangent current to $T$ along $V$ associated with $\lambda_n$. 
Denote by $S$ this tangent current and write $S=\pi^*(\nu)$ as in Section \ref{section_currents}. 
By Theorem \ref{th_tangent}, it is enough to show that $\nu=0$. 
We will see later that the choice of $\lambda_n$ is very important for the techniques we use below.
Recall that by Proposition \ref{prop_tangent_local}, we can compute $S$ and $\nu$ using local coordinates. 

\begin{lemma}
We have $\nu=0$ outside the singular points of $\Fc$ in $V$. 
\end{lemma}
\proof
Let $p\in V$ be a regular point for $\Fc$. We show that $\nu=0$ in a neighbourhood of $p$. 
Choose local coordinates $z=(z_1,z_2)$ centered at $p$, $|z_1|<3$, $|z_2|<3$, on which $\Fc$ is given by the horizontal lines  $\{z_2=\alpha\}$ with $\alpha\in\C$. In particular, $V$ is given there by $z_2=0$. In this chart, the current $T$ can be written as
$$T=\int h_\alpha [z_2=\alpha] d\mu(\alpha),$$
where $\mu$ is a positive measure on the disc $\D(0,3)$ of $\C$ and $h_\alpha$ is a positive harmonic function in $z_1\in \D(0,3)$. 
Dividing  $h_\alpha$ by $h_\alpha(0)$ and multiplying $\mu$ by $h_\alpha(0)$, we can assume that $h_\alpha(0)=1$ for $\mu$-almost every $\alpha\in\A$. It follows from Harnack's inequality that on $\D$ the harmonic function $h_\alpha$ is bounded below and above by strictly positive constants independent of $\alpha$. We also deduce that $\mu$ has finite mass on $\D$. 
Since $T$ has no mass on $V$, we can assume that $\mu$ has no mass at 0.

We can now apply Proposition \ref{prop_tangent_local} to the considered coordinates. We obtain that 
$$T_\lambda=\int h_{\alpha/\lambda} [z_2=\alpha] d\mu_\lambda(\alpha),$$
where $\mu_\lambda$ is the direct image of $\mu$ by the dilation $\alpha\mapsto\lambda\alpha$. 
 The fact that $\mu$ has no mass at 0 implies that the mass of $\mu_\lambda$ in any compact set tends to 0 as $\lambda\to\infty$. 
Finally, since $h_\alpha$ are uniformly bounded on $\alpha$, we conclude that $T_\lambda$ converges to 0. The lemma follows. 
\endproof

The last lemma shows that $\nu$ is a finite combination of Dirac masses at singular points of $\Fc$ in $V$. 
In what follows, consider a singular point $p\in V$ of $\Fc$. It remains to prove that $\nu(\{p\})=0$. 
We use the notation introduced in the previous sections. 
Fix a smooth function $\vartheta(z_2)$ with compact support in the annulus $\{e^{-1}\leq |z_2|\leq 1\}$
such that $\int \vartheta(z_2)idz_2\wedge d\overline z_2=1$. It is enough to check that the mass of the measure 
$S\wedge \vartheta(z_2) idz_2\wedge d\overline z_2$ in a neighbourhood of $\{0\}\times\C$ is zero. 
Indeed, the properties of tangent currents developed in Theorem \ref{th_tangent} imply that this mass is also the mass of $\nu$ at the singular point $p$. 
By definition of $S$, it suffices to show that the mass of $T\wedge |\lambda_n|^2 \vartheta(\lambda_n z_2) idz_2\wedge d\overline z_2$ on $\D\times \C$ tends to 0 as $n\to\infty$. 

Define for $s\geq 1$
$$G_\alpha(s):=\int_{u+is\in \S} H_\alpha(u+is)du.$$
We have the following lemma.

\begin{lemma} \label{lemma_T_alpha}
For $\lambda=e^s$ with $s\geq 1$, the mass of $T_\alpha\wedge |\lambda|^2 \vartheta(\lambda z_2) idz_2\wedge d\overline z_2$ on $\overline \D\times \C$ is bounded above by a constant times $G_\alpha(s)$.
\end{lemma}
\proof
Note that the support of $\vartheta(z_1,\lambda z_2)$ is contained in $ \{e^{-s-1}\leq |z_2|\leq e^{-s}\}$.
Since $\vartheta$ is bounded, the mass considered in the lemma is bounded by a constant times the mass of the positive measure $T_\alpha\wedge \lambda^2 idz_2\wedge d\overline z_2$ on $\overline \D\times \{e^{-s-1}\leq |z_2|\leq e^{-s}\}$. We will only consider $z_2$ such that $e^{-s-1}\leq|z_2|\leq e^{-s}$. Recall that $T_\alpha=h_\alpha[L_\alpha]$. So the last mass is equal to
$$\int_{e^{-s-1}\leq |z_2|\leq e^{-s}} \Big[\sum_{(z_1,z_2)\in L_\alpha} h_\alpha(z_1,z_2)\Big] e^{2s} idz_2\wedge d\overline z_2.$$
Since $e^{2s} idz_2\wedge d\overline z_2$ is a measure of bounded mass in $\{e^{-s-1}\leq |z_2|\leq e^{-s}\}$, we only need to bound the sum in the previous brackets by a constant times $G_\alpha(s)$.  

For each $z_2$, we consider the sequence of points $(z_1,z_2)$ in $L_\alpha$.
Write $z_2=e^{-v+iu_0+i\log|\alpha|/b}$ with  $v:=-\log|z_2|$ and $u_0\in \R$ such that 
$$-{av \over b} \leq   u_0 <  -{av \over b}+2\pi.$$
We have $s\leq v\leq s+1$ and hence $-as/b-|a|/b\leq u_0< -as/b+|a|/b+2\pi$. So both $v$ and $u_0$  vary in intervals of bounded length.

Consider the set
$$I  :=   \big\{\zeta_k:=u_k+iv \quad \text{with}\quad u_k=u_0+2k\pi,\quad k=0,1, 2, \ldots \big\}.$$
This is an arithmetic progression on a half-line in $\overline\S$. 
We have 
$$\sum_{(z_1,z_2)\in L_\alpha} h_\alpha(z_1,z_2)=\sum_{k\geq 0} H_\alpha(u_k+iv).$$
Fix a constant $A>0$ large enough. Recall that $H_\alpha$ is a positive harmonic function in the sector $\S'$ which contains the sector $\overline \S$. By Harnack's inequality, $H_\alpha(u_k+iv)$ for $s\leq v\leq s+1$
  is bounded above by a  positive constant times $H_\alpha(u+is)$ for all  $u$ in $[u_k-A,u_k+A]$ with $u+is\in\S$. We then deduce that the last infinite sum is bounded by a constant times $G_\alpha(s)$. The lemma follows.
\endproof

Denote by $\Lambda$ the half-line $\S\cap\{v=1\}$.
Let $\zeta^*=u^*+i$ be the intersection point of $\Lambda$ with the boundary of $\S$. Define also $\rho:=|\zeta^*|^\gamma=-\Phi(\zeta^*)$ which is a positive number depending only on $\gamma$.
We will use Lemma \ref{lemma_s_n_bis} for this value of $\rho$.  
The image of $\Lambda$ by $\Phi$ is denoted by $\widetilde \Lambda$. It is an unbounded curve in the half-plane $\H$ starting from the point $-\rho$. We will transform the integrals to be estimated into integrals on $\Lambda$ or $\widetilde\Lambda$. For this purpose, we  
use the variable $r:=u/s-u^*$ with $r\geq 0$. 
Since we will consider  $\zeta$ in the line $v=s$, we write 
 $$\zeta=u+is, \qquad Z=U+iV=\Phi(\zeta)=\zeta^\gamma=(u+is)^\gamma$$
and 
$$Z'=s^{-\gamma} Z, \qquad U'=s^{-\gamma} U, \qquad V'=s^{-\gamma} V.$$
The variables $Z,U,V$ depend on $s,r$ and the variables $Z',U',V'$ can be seen as functions in $r$. The point $Z'=U'+iV'$ belongs to $\widetilde\Lambda$. 
We obtain with these notations and using Lemma \ref{lemma_poisson} that
$$G_\alpha(s)=\int_0^\infty \widetilde H_\alpha(U+iV) sdr={1\over \pi}\int_0^\infty  \int_{x\in\R} \widetilde H_\alpha(x){V\over V^2+(U-x)^2} s dx dr.$$
Write also $x=s^\gamma x'$. We have
$$G_\alpha(s)={1\over \pi}\int_0^\infty  \int_{x'\in\R} \widetilde H_\alpha(x){V'\over V'^2+(U'-x')^2} s dx' dr.$$

Consider a constant $N\geq 1$ large enough whose value will be specified later. We split the last integral into two parts and write
$$G_\alpha(s)=G_{\alpha,N}^{(1)}(s) + G_{\alpha,N}^{(2)}(s),$$
where
$$G_{\alpha,N}^{(1)}(s) ={1\over \pi}\int_0^\infty  \int_{|x'|\geq N} \widetilde H_\alpha(x){V'\over V'^2+(U'-x')^2} s dx' dr$$
and
$$G_{\alpha,N}^{(2)}(s) ={1\over \pi}\int_0^\infty  \int_{|x'|\leq N} \widetilde H_\alpha(x){V'\over V'^2+(U'-x')^2} s dx' dr.$$
The aim is to bound these integrals. We first give some properties of $U'$ and $V'$. Recall  that $Z'=U'+iV'$ describes the curve $\widetilde \Lambda$ when $r$ varies in $[0,\infty)$. 

\begin{lemma} \label{lemma_Phi_Lambda_0}
We have that  $U'=-\rho+O(r)$ and $V'=\beta r+O(r^2)$ for some constant $\beta>0$ when $r\to 0$. 
Moreover, given a constant $N>0$, we have for $0\leq r\leq N$
$$\dist(x',U'+iV')^2\geq c_N \Big[r^2+\dist(x', -\rho)^2\Big],$$
where $\dist$ denotes the standard distance and $c_N>0$ is a constant independent of $x'$.
\end{lemma}
\proof
When $r\to 0$, we have 
$$Z'=U'+iV'=\Phi(\zeta^*+r) =-\rho+c_1r+c_2r^2+O(r^3)$$
for some constants $c_1$ and $c_2$. So clearly $U'=-\rho+O(r)$. 
Moreover, observe that $\Phi$ is conformal near $\zeta^*$. Since the boundary of $\S$ is sent by $\Phi$ to $\R$ and the half-line $\Lambda$ intersects the boundary of $\S$ transversally at $\zeta^*$, we deduce that the curve $\widetilde\Lambda=\Phi(\Lambda)$ intersects $\R$ transversally at $-\rho$. It follows easily that $V'=\beta r+O(r^2)$ for some constant $\beta>0$.

For the second assertion, let $0<\theta<\pi$ be the angle at $-\rho$ of the triangle of vertices $-\rho,U'+iV'$ and $x'$.  
By the law of cosines, we have 
$$\dist(x',U'+iV')^2 = \dist(x',-\rho)^2+\dist (-\rho,U'+iV')^2 -2\cos \theta \dist(x',-\rho)\dist(-\rho,U'+iV').$$
The right hand side is larger or equal to 
$$(1-\cos\theta)\Big[\dist(x',-\rho)^2+\dist (-\rho,U'+iV')^2\Big].$$
When $0\leq r\leq N$, the angle $\theta$ is bounded below by a strictly positive constant. Therefore, $1-\cos\theta$ is bounded below by a strictly positive constant. The first assertion implies that $\dist (-\rho,U'+iV')\gtrsim r$. The result follows.
\endproof

\begin{lemma} \label{lemma_Phi_Lambda_infty}
When $r\to\infty$, we have
$$U'=r^\gamma+O(r^{\gamma-1}) \qquad \text{and} \qquad V'=\gamma r^{\gamma-1}+ O(r^{\gamma-2}).$$
\end{lemma}
\proof
Observe that $\Phi(\zeta)=\zeta^\gamma$ and  $Z'=\Phi(u^*+r+i)=r^\gamma\Phi(1+u^*/r+i/r)$. Therefore, when $r\to\infty$, we have that $Z'=r^\gamma+O(r^{\gamma-1})$. Since the argument of $u^*+r+i$ is $1/r+O(1/r^2)$, the argument of $Z'$ is $\gamma/r+O(1/r^2)$.  It follows that $rV'/U'=\gamma+O(1/r)$  when $r\to\infty$. We then easily deduce desired estimates for $U'$ and $V'$ from the above estimate for $Z'$.
\endproof

\begin{lemma} \label{lemma_dr}
There is a constant $c>0$ such that for $|x'|\geq 2\rho$
$$ \int_0^\infty {V'\over V'^2+(U'-x')^2} dr \leq c |x'|^{-1+1/\gamma}.$$
\end{lemma}
\proof
We split this integral into two parts : an integral on $[0,1]$ and another on $[1,\infty)$. Consider first the case where $0\leq r\leq 1$. 
Recall that $U'+iV'$ is in the curve $\widetilde\Lambda=\Phi(\Lambda)$. Since $|x'|\geq 2\rho$, 
we have $|x'+\rho| \gtrsim |x'|$. Using the  last estimate in Lemma \ref{lemma_Phi_Lambda_0}, we get 
$$ \int_0^1 {V'\over V'^2+(U'-x')^2} dr \lesssim  \int_0^1 {dr \over |x'+\rho|^2} \lesssim \int_0^1 {dr \over |x'|^2}= |x'|^{-2}\lesssim |x'|^{-1+1/\gamma}.$$

It remains to prove a similar estimate for the integral on $[1,\infty)$. By Lemma \ref{lemma_Phi_Lambda_infty}, we only need to check that 
$$ \int_1^\infty {r^{\gamma-1} dr\over r^{2\gamma-2}+(U'-x')^2} \leq c |x'|^{-1+1/\gamma}$$
for some constant $c>0$. It is not difficult to see using again Lemma \ref{lemma_Phi_Lambda_infty} that the last integral is bounded when $x'$ is bounded. Therefore, we only need to consider $|x'|$ large enough.

Fix a constant $A>1$ large enough depending only on $\gamma$. We can assume that  $|x'|\gg A^{3\gamma}$. By Lemma \ref{lemma_Phi_Lambda_infty}, we have $U'\geq r^\gamma -Ar^{\gamma-1}$. Therefore, if $r\geq |x'|^{1/\gamma}+A^3$, then 
$|x'|\leq (r-A^3)^\gamma\leq (r-A^2)^\gamma$.
Using the mean value theorem for the function $t\mapsto t^\gamma$ between $t=r$ and $t=r-A^2$,  we obtain
$${1\over 2}\big(r^\gamma-|x'|\big)\geq {1\over 2}\big(r^\gamma-(r-A^2)^\gamma\big)\geq {1\over 2}\gamma A^2(r-A^2)^{\gamma-1}=
Ar^{\gamma-1}{\gamma\over 2} A(1-A^2/r)^{\gamma-1}\geq Ar^{\gamma-1}$$
since $A^2/r\ll 1/A \ll 1$, and hence
$$U'-x'\geq (r^\gamma-Ar^{\gamma-1})-|x'|=(r^\gamma-|x'|)-Ar^{\gamma-1}\geq {1\over 2} (r^\gamma-|x'|).$$
It follows that 
\begin{eqnarray*}
\int_{r\geq |x'|^{1/\gamma}+A^3} {r^{\gamma-1} dr\over r^{2\gamma-2}+(U'-x')^2}  &\lesssim &  
  \int_{r\geq |x'|^{1/\gamma}+A^3} {dr^\gamma \over (r^\gamma-|x'|)^2} \\
 &\lesssim & {1\over (|x'|^{1/\gamma}+A^3)^\gamma-|x'|}\\ 
 & \lesssim &  |x'|^{-1+1/\gamma},
 \end{eqnarray*}
where we use again the mean value theorem for the function $t\mapsto t^\gamma$.

Similarly, we have $|U'|\leq r^\gamma+Ar^{\gamma-1}$. If $1\leq r\leq |x'|^{1/\gamma}-A^3$,  
using the mean value theorem always for the function $t\mapsto t^\gamma$ but between 
$|x'|^{1/\gamma}$ and $|x'|^{1/\gamma}-A^2$, we get
$$|x'-U'|\geq |x'|-(r^\gamma+Ar^{\gamma-1})= (|x'|-r^\gamma)-Ar^{\gamma-1}\geq {1\over 2} (|x'|-r^\gamma)$$
and we deduce in the same way that
$$\int_{1\leq r\leq |x'|^{1/\gamma}-A^3} {r^{\gamma-1} dr\over r^{2\gamma-2}+(U'-x')^2}   \lesssim   |x'|^{-1+1/\gamma}.$$

On the other hand, we have (recall that we only consider $|x'|$ large enough)
$$\int_{|x'|^{1/\gamma}-A^3\leq r\leq |x'|^{1/\gamma}+A^3} {r^{\gamma-1} dr\over r^{2\gamma-2}+(U'-x')^2}  \lesssim 
\int_{|x'|^{1/\gamma}-A^3\leq r\leq |x'|^{1/\gamma}+A^3} {r^{\gamma-1}dr \over r^{2\gamma-2}} 
  \lesssim   |x'|^{-1+1/\gamma}$$
since $r\simeq |x'|^{1/\gamma}$ on the interval where we take the integral and this interval has a bounded length.  The last estimate and the two previous ones imply the lemma.
\endproof

\begin{lemma} \label{lemma_G_1}
Let  $\epsilon>0$ be an arbitrary fixed constant. If the constant $N$ is large enough, then we have for every $s\geq 1$
$$\int_{\alpha\in\A} G_{\alpha,N}^{(1)}(s)d\mu(\alpha)\leq \epsilon.$$ 
\end{lemma}
\proof
Using Lemma \ref{lemma_dr} and that  $s\geq 1$, we can bound $G_{\alpha,N}^{(1)}(s)$ by a constant times
$$\int_{|x'|\geq N} \widetilde H_\alpha(x) |x'|^{-1+1/\gamma} sdx'
= \int_{|x|\geq Ns^\gamma} \widetilde H_\alpha(x) |x|^{-1+1/\gamma} dx
\leq \int_{|x|\geq N} \widetilde H_\alpha(x) |x|^{-1+1/\gamma} dx.$$
Denote by $\widetilde G_{\alpha,N}^{(1)}$ the last integral which is independent of $s$. 
By Lemma \ref{lemma_poisson_est}, $\widetilde G_{\alpha,N}^{(1)}$ is bounded by a constant independent of $s,\alpha$ and $N$. Moreover, for $\mu$-almost every $\alpha\in\A$ we have
$$\lim_{N\to\infty} \widetilde G_{\alpha,N}^{(1)}=0.$$
By Lebesgue's convergence theorem, we obtain
$$\lim_{N\to\infty} \int_{\alpha\in\A}  \widetilde G_{\alpha,N}^{(1)}d\mu(\alpha)=0.$$
The lemma follows.
\endproof

\begin{lemma} \label{lemma_G_2}
Let $N\geq 2\rho$ be a fixed constant large enough. Then we have 
$$\lim_{n\to\infty} \int_{\alpha\in \A} G_{\alpha,N}^{(2)}(s_n) d\mu(\alpha)= 0.$$
\end{lemma}
\proof
We only consider large $n$. From the definition of $G_{\alpha,N}^{(2)}$, we get
\begin{eqnarray*}
G_{\alpha,N}^{(2)}(s_n) & = &  {1\over \pi}s_n\int_0^\infty \int_{|x'|\leq N} \widetilde H_\alpha(x) {V'\over V'^2+(U'-x')^2} dx'dr \\
& = &  {1\over \pi}s_n^{1-\gamma}\int_0^\infty \int_{|x|\leq Ns_n^\gamma} \widetilde H_\alpha(x) {V'\over V'^2+(U'-x')^2} dxdr.
\end{eqnarray*}
Therefore, by Lemma \ref{lemma_s_n_bis}, we only need to check for $|x'|\leq N$  that 
$$\int_0^\infty {V'\over V'^2+(U'-x')^2} dr \lesssim \log {6Ns^\gamma\over |x+\rho s^\gamma|}$$
or equivalently
$$\int_0^\infty {V'\over V'^2+(U'-x')^2} dr \lesssim \log {6N\over |x'+\rho|}\cdot$$

Observe that the right hand side of the last inequality is larger than 1 because $|x'|\leq N$ and $N\geq 2\rho$. Since $x'$ is bounded, Lemma \ref{lemma_Phi_Lambda_infty} implies that 
$$\int_1^\infty {V'\over V'^2+(U'-x')^2} dr \lesssim \int_1^\infty {r^{\gamma-1}dr\over r^{2\gamma}}= \const \lesssim \log {6N\over |x'+\rho|}\cdot$$
 It remains to check that 
$$\int_0^1 {V'\over V'^2+(U'-x')^2} dr \lesssim \log {6N\over |x'+\rho|}\cdot$$
Using the last assertion in Lemma \ref{lemma_Phi_Lambda_0}, we have
\begin{eqnarray*}
\int_0^1 {V'\over V'^2+(U'-x')^2} dr & \lesssim & \int_0^1 {rdr\over r^2+|x'+\rho|^2} \ = \ {1\over 2} \log{1+|x'+\rho|^2 \over|x'+\rho|^2}\\
& = & \log {\sqrt{1+|x'+\rho|^2}\over |x'+\rho|} \  \lesssim \  \log {6N\over |x'+\rho|}
\end{eqnarray*}
because $|x'+\rho|\leq 2N$. 
This ends the proof of the lemma.
\endproof

\noindent
{\bf End of the proof of Theorem \ref{th_main_bis}.} 
Recall that $\lambda_n=e^{s_n}$. 
We have to show that the mass of $T\wedge |\lambda_n|^2\vartheta(\lambda_n z_2) idz_2 \wedge d\overline z_2$ on $\overline\D\times \C$ tends to 0 as $n\to\infty$. 
By Lemma \ref{lemma_T_alpha}, this mass  is bounded by a constant times 
$$\int_{\alpha\in \A} G_\alpha(s_n)  d\mu(\alpha)=\int_{\alpha\in \A}  G_{\alpha,N}^{(1)}(s_n)  d\mu(\alpha)+\int_{\alpha\in \A}  G_{\alpha,N}^{(2)}(s_n) d\mu(\alpha).$$
Let $\epsilon>0$ be an arbitrary constant. Fix also an $N$ large enough. Lemmas \ref{lemma_G_1} and \ref{lemma_G_2} imply that the upper limit of the last sum is bounded above by $\epsilon$. This ends the proof in the case where $V$ is smooth.
In general, the only possible singularities of $V$ are the singularities of $\Fc$ near which $V$ is the union of the two separatrices. 
Our study can be extended to this case without difficulty. We present here an alternative argument. 

Consider the blow-up $\pi:\widehat X\to X$ of $X$ at the singularities of $V$. Let $E_1,\ldots, E_m$ denote the exceptional fibers of $\pi$ and $\widehat V$ the strict transform of $V$. The foliation $\Fc$ can be lifted to a foliation $\widehat\Fc$. It is not difficult to check that $\widehat V$ and $E_i$ are smooth irreducible invariant curves for $\widehat\Fc$ and all singularities of $\widehat\Fc$ on these curves are hyperbolic, see \cite[p.6]{Brunella}.
Let $\widehat T$ denote the pullback of $T$ by $\pi$ to $\widehat X\setminus\cup E_i$ which is a harmonic current of $\widehat \Fc$ without mass on $\widehat V$ and $E_i$. By the first case considered above, we get $\{\widehat T\}\smallsmile \{\widehat V\}=0$ and $\{\widehat T\}\smallsmile \{E_i\}=0$. It follows that $\{\widehat T\}\smallsmile \pi^*(\{V\})=0$ because $\pi^*(\{V\})$ is a combination of $\{\widehat V\}$ and $\{E_i\}$. 
Thus, $\{T\}\smallsmile \{V\}=0$ because
$\pi_*(\widehat T)=T$.
\hfill $\square$

\end{document}